\documentclass[10pt]{amsart}
\usepackage[english]{babel}
\usepackage{amsmath}
\usepackage{amssymb}
\usepackage{amsthm}
\usepackage{libertine}
\usepackage[all]{xy}
\usepackage{tikz-cd}
\usepackage{tikz}
\usepackage{float}
\newtheorem{thm}{Theorem}[section]
\newtheorem{cor}[thm]{Corollary}
\newtheorem{lem}[thm]{Lemma}
\newtheorem{prop}[thm]{Proposition}

\theoremstyle{definition}
\newtheorem{defn}[thm]{Definition}
\theoremstyle{remark}
\newtheorem{rem}[thm]{Remark}

\numberwithin{equation}{section}
\newcommand{\norm}[1]{\left\Vert#1\right\Vert}
\newcommand{\scal}[1]{\left<#1\right>}
\newcommand{\Hq}{\mathbb H}
\newcommand{\Sq}{\mathbb S}
\newcommand{\R}{\mathbb{R}}      

\newcommand{\C}{\mathbb{C}}

\title[Generalized Appell polynomials and Fueter-Bargmann transforms in the polyanalytic setting]{Generalized Appell polynomials and Fueter-Bargmann transforms in the polyanalytic setting}

\begin{document}
\date{}
\author{Antonino De Martino, \, Kamal Diki }
\maketitle
\begin{abstract}
This paper deals with some special integral transforms in the setting of quaternionic valued slice polyanalytic functions. In particular, using the polyanalytic Fueter mappings it is possible to construct a new family of polynomials which are called the generalized Appell polynomials. Furthermore, the range of the polyanalytic Fueter mappings on two different polyanalytic Fock spaces is characterized. Finally, we study the polyanalytic Fueter-Bargmann transforms.
\end{abstract}


\section{Introduction}
This paper proposes a generalization to the polyanalytic setting of both Appell polynomials and the Bargmann-Fock-Fueter transform, studied in \cite{CMF2017, CMF2011} and \cite{DKS2019} respectively. A fundamental tool to perform the computations is the polyanalytic Fueter map, introduced in \cite{ADS}. One of the main differences with respect to the classic Fueter theorem, see \cite{F}, is that in the polyanalytic setting there are two Fueter mappings. The first map, denoted by $ \mathcal{C}_{n+1}$, maps slice polyanalytic functions of order $n+1$ to Fueter regular functions of order $n+1$. The second one, denoted by $ \tau_{n+1}$, maps slice polyanalytic functions of order $n+1$ to Fueter regular ones. We note that there exists a relation between these two polyanalytic Fueter maps that can be expressed in terms of a suitable power of the Cauchy-Fueter operator, see \cite{ADS}.
\\In Section 3 we introduce a new family of polynomials in the quaternionic setting
\begin{equation}
\label{pp}
\mathcal{M}_{k,s}(q,\overline{q}):=x_0^kQ_s(q,\overline{q}),\quad k=0,...,n, \quad s\geq 0,
\end{equation}
where $ Q_s(q,\overline{q})=\sum_{j=0}^{s} \frac{2(s-j+1)}{(s+1)(s+2)}q^{s-j}\overline{q}^j$ and $x_0$ is the real part of the quaternion $q \in \mathbb{H}$. 
These polynomials are obtained by applying the polyanalytic Fueter map $\mathcal{C}_{n+1}$ to a slice polyanalytic function of order $n+1$, based on the series expansion theorem. One of the main properties of these polynomials is the following 
$$\overline{\mathcal{D}}^{n+1}(\mathcal{M}_{n,s}(q,\overline{q}))=\displaystyle\sum_{j=1}^{n+1} 2^j {n+1\choose j}\frac{n!s!}{(j-1)!(s-j)!}\mathcal{M}_{j-1,s-j}(q,\overline{q}), \quad s \geq n+1, $$
where $\overline{\mathcal{D}}$ is the hypercomplex derivate.
The previous formula, when $n=0$, leads to the classical Appell property for the Clifford Appell polynomials in the quaternionic setting
$$ \overline{\mathcal{D}} (\mathcal{Q}_s(q, \bar{q}))=2s \mathcal{Q}_{s-1}(q, \bar{q}),$$
see \cite{CMF2017,CMF2011, DKS2019}.
However, when we apply the second polyanalytic Fueter map $ \tau_{n+1}$, we do not have any suggestions about new family of polynomials. This is due to the fact that the second polyanalytic Fueter map has a range included in the space of Fueter regular functions. In \cite{DDG1}, by using some results proved in \cite{DDG}, is provided a generalization of the polynomials defined in formula \eqref{pp} to the Clifford Algebra in the odd dimensions.
\\ In Section 4 we investigate the polyanalytic Fueter-Bargmann transforms. Firstly, we study a characterization of the quaternionic true polyanalytic Fock space $ \mathcal{F}_n^T(\mathbb{H})$, see \cite{DMD2}, when the Fueter map $ \mathcal{C}_{n+1}$ is applied. After, we define the $ \mathcal{C}$- polyanalytic Fueter-Bargmann transform as
$$ B^{n+1}_{\mathcal{C}}:= \mathcal{C}_{n+1}\circ B^{n+1},$$
where $B^{n+1}$ is the quaternionic true polyanalytic Bargmann transform defined in \cite{DMD2}. By using the fact that it is possible to write the kernel of the quaternionic true polyanalytic Bargmann transform as a generating function of the quaternionic Hermite polynomials (see \cite{TT}), an explicit integral formula for $B^{n+1}_{\mathcal{C}}$ is showed. Besides, some important properties for the $ \mathcal{C}$- polyanalytic Fueter-Bargmann transform, such as the unitary and isometric properties are showed. In the second part of Section 4 the polyanalytic Fueter map $ \tau_{n+1}$ is applied. In this case, most of the results are obtained by using the connection between the maps $\mathcal{C}_{n+1}$ and $\tau_{n+1}$.
\\ Finally, in Section 5 the results proved in the previous section are discussed for the quaternionic polyanalytic Fock space $ \mathcal{F}^{N+1}_{Slice}(\mathbb{H})$, studied in \cite{ADS2019, DMD2}.

\section{Preliminary results}
We revise different notions and results about quaternions and related function theories. The non-commutative field of quaternions is defined to be
$$\Hq=\lbrace{q=x_0+x_1i+x_2j+x_3k : \ x_0,x_1,x_2,x_3\in\R}\rbrace,$$ where the imaginary units satisfy the multiplication rules $$i^2=j^2=k^2=-1\quad \text{and}\quad ij=-ji=k, jk=-kj=i, ki=-ik=j.$$
On $\Hq$ the conjugate and the modulus of $q$ are defined respectively by
$$\overline{q}=x_0-\vec{q}\,,  \quad \vec{q}\,=x_1i+x_2j+x_3k$$
and $$\vert{q}\vert=\sqrt{q\overline{q}}=\sqrt{x_0^2+x_1^2+x_2^2+x_3^2}.$$
We note that the quaternionic conjugation satisfy the property $\overline{ pq }= \overline{q}\, \overline{p}$ for any $p,q\in \Hq$.
Moreover, the unit sphere $$\lbrace{q=x_1i+x_2j+x_3k : \text{ } x_1^2+x_2^2+x_3^2=1}\rbrace$$ coincides with the set of all  imaginary units given by $$\mathbb{S}=\lbrace{q\in{\Hq} : q^2=-1}\rbrace.$$
Any quaternion $q\in \Hq\setminus \R$ can be written in a unique way as $q=x+I y$ for some real numbers $x$ and $y>0$, and imaginary unit $I\in \mathbb{S}$. For every given $I\in{\mathbb{S}}$, we define $\C_I = \mathbb{R}+\mathbb{R}I.$
It is isomorphic to the complex plane $\C$ so that it can be considered as a complex plane in $\Hq$ passing through $0$, $1$ and $I$. Their union is the whole space of quaternions $$\Hq=\underset{I\in{\mathbb{S}}}{\cup}\C_I =\underset{I\in{\mathbb{S}}}{\cup}(\mathbb{R}+\mathbb{R}I).$$
Let $\mathbb{B}$ denotes the quaternionic unit ball and $\mathbb{B}_I$ its intersection with the complex plane $\C_I$ for a given $I\in\mathbb{S}$. Then, we recall
\begin{defn} Let $U\subset \Hq$ be an open set and let $f:U\longrightarrow \Hq$ be a function of class $\mathcal{C}^1$. We say that $f$ is (left) Fueter regular (or regular for short) on $U$ if $$\mathcal{D} f(q):=\displaystyle\left(\frac{\partial}{\partial x_0}+i\frac{\partial}{\partial x_1}+j\frac{\partial}{\partial x_2}+k\frac{\partial}{\partial x_3}\right)f(q)=0, \quad \forall q\in U.$$
\\
The right quaternionic vector space of Fueter regular functions will be denoted by $\mathcal{FR}(U)$.
\end{defn}
We recall that a slice domain (or s-domain) $U$ on $\mathbb{H}$ is a domain that intersect the real line and for which $U_I=U\cap \mathbb{C}_I$ is a domain of $\mathbb{C}_I$ for every $I\in\mathbb{S}$. Moreover, if for every $q=x+Iy\in{\Omega}$, the whole sphere $x+y\mathbb{S}:=\lbrace{x+Jy; \, J\in{\mathbb{S}}}\rbrace$
is contained in $\Omega$, we say that  $\Omega$ is an axially symmetric slice domain. 
\begin{defn}\label{slicefunction}
Let $U$ be an axially symmetric set in $\mathbb H$. A quaternionic valued function $f:U\subset \Hq\longrightarrow \Hq$ is called a slice function if it is of the form $f(x+yI)=\alpha(x,y)+I\beta(x,y)$, where $\alpha(x,y)$ and $\beta(x,y)$ are quaternionic-valued functions such that $\alpha(x,-y)=\alpha(x,y)$, $\beta(x,-y)=-\beta(x,y).$
\end{defn}
In \cite{GS2007} the theory of slice hyperholomorphic functions was introduced based on the following notion:
	
\begin{defn}
Let $f: \Omega \longrightarrow \Hq$ be a  $\mathcal{C}^1$ slice function on a given domain $\Omega\subset \Hq$. Then, $f$ is said to be (left) slice hyperholomorphic function if, for every $I\in \Sq$, the restriction $f_I$ to $\C_{I}=\R+I\R$, with variable $q=x+Iy$, is holomorphic on $\Omega_I := \Omega \cap \C_I$, that is it has continuous partial derivatives with respect to $x$ and $y$ and the function
$\overline{\partial_I} f : \Omega_I \longrightarrow \Hq$ defined by
$$
\overline{\partial_I} f(x+Iy):=
\dfrac{1}{2}\left(\frac{\partial }{\partial x}+I\frac{\partial }{\partial y}\right)f_I(x+yI)
$$
vanishes identically on $\Omega_I$. The set of such kind of functions will be denoted by $ \mathcal{SR}(\Omega)$.
\end{defn}
\begin{defn}
The slice derivative of a slice regular function $f$ is defined as:
\begin{equation*}
\partial_S(f)(q):=
\left\{
\begin{array}{rl}
\partial_I(f)(q)& \text{if } q=x+Iy, y\neq 0\\
\frac{\partial}{\partial{x}}(f)(x) & \text{if } q=x \text{ is real}.
\end{array}
\right.
\end{equation*}
\end{defn}
The right quaternion vector space of slice hyperholomorphic functions is endowed with the natural topology of uniform convergence on compact sets. In particular, we recall the following series expansion from \cite[Cor. 4.2.3]{CSS1}.
\begin{thm}
Let $f$ be a slice hyperholomorphic function on an axially symmetirc s-domain $U$. Then for any real point $p_0$ in $U$, the function $f$ can be represented by power series
$$f(q)=\sum_{m=0}^{+\infty} (q-p_0)^ma_m$$
on the ball $B(p_0,R)=\{q\in\Hq; |q-p_0|<R\}$ where $R = R_{p_0}$ is the largest positive real number such that
$B(p_0,R)$ is contained in $U$.
\end{thm}
These two quaternionic function theories can be related using a fundamental result in quaternionic analysis which is the Fueter mapping theorem. We recall briefly this result here
\begin{thm}[Fueter mapping theorem, \cite{F}]
Let $U$ be an axially symmetric set in $\mathbb H$ and let $f:U\subset \Hq\longrightarrow \Hq$ be a slice regular function of the form $f(x+yI)=\alpha(x,y)+I\beta(x,y)$ as in Definition \ref{slicefunction} and satisfying the Cauchy-Riemann system. Then, the function $$\overset{\sim}f(x_0,\vert{\vec{q}}\vert)=\displaystyle \Delta\left(\alpha(x_0,\vert{\vec{q}}\vert)+\frac{\vec{q}}{\vert{\vec{q}}\vert}\beta(x_0,\vert{\vec{q}}\vert)\right)$$
is Fueter regular.
\end{thm}
In the late of 1950s, Sce extended this theorem to the Clifford setting in the case of odd dimensions see \cite{CSS3, S}, Qian proved in \cite{Q} that the theorem of Sce holds in the case of even dimensions. We refer to \cite{ColomboSabadiniSommen2010, Q1} for several extensions.
	
\begin{rem}
We denote the Fueter mapping by $$\tau:\mathcal{SR}(U)\rightarrow\mathcal{FR}(U), \text{ } f\longmapsto \tau(f)=\overset{\sim}f.$$
\end{rem}
We note that one can extend the classical theory of holomorphic functions to higher order using the notion of polyanalytic functions, see \cite{B}. Moreover, in the last years the notion of slice hyperholomorphic functions was extended also to this polyanalytic setting, see \cite{ACDS1, ADS2019, ADS, BEA}. We briefly recall this notion here and related results that will be needed in the sequel
\begin{defn}
Let $\Omega$ be an axially symmetric open set in $\Hq$ and let $f:\Omega\longrightarrow \Hq$ a slice function of class $\mathcal{C}^{n+1}$. For each $I\in\Sq$, let $\Omega_I=\Omega\cap\C_I$ and let $f_I=f_{\vert_{\Omega_I}}$ be the restriction of $f$ to $\Omega_I$. The restriction $f_I$ is called (left) polyanalytic of order $n+1$ if it satisfies on $\Omega_I$ the equation $$
\overline{\partial_I}^{n+1} f(x+Iy):=
\frac{1}{2^{n+1}}\left(\frac{\partial }{\partial x}+I\frac{\partial }{\partial y}\right)^{n+1}f_I(x+Iy)=0.
$$
The function $f$ is called left slice polyanalytic of order $n+1$, if for all $I\in\Sq$, $f_I$ is left polyanalytic of order $n+1$ on $\Omega_I$. The right quaternionic vector space of slice polyanalytic functions of order $n+1$ will be denoted by $\mathcal{SP}_{n+1}(U)$.
\end{defn}
Note that slice regular functions are a special case of the definition of slice polyanalytic functions with $n=0$. Several results of these functions were studied and extended. In particular, we note that the following decomposition holds true
\begin{prop}[polyanalytic-decomposition]\label{carac3}
A slice function $f:\Omega\longrightarrow \Hq$ defined on an axially symmetric slice domain is slice polyanalytic of order $n+1$ if and only if there exist $f_0,...,f_{n}$ some unique slice hyperholomorphic functions on $\Omega$ such that we have the following decomposition: $$
f(q):=\sum_{k=0}^{n}\overline{q}^kf_k(q); \textbf{ }\forall q\in\Omega.
$$
\end{prop}
	
We revise also the polyanalytic Fueter regular functions (see \cite{B1976, DB1978}).
\begin{defn} Let $U\subset \Hq$ be an open set and let $f:U\longrightarrow \Hq$ be a function of class $\mathcal{C}^{n+1}$. We say that $f$ is (left) polyanalytic Fueter regular (or polyanalytic-regular for short) of order $n+1$ on $U$ if $$\mathcal{D}^{n+1} f(q):=\displaystyle\left(\frac{\partial}{\partial x_0}+i\frac{\partial}{\partial x_1}+j\frac{\partial}{\partial x_2}+k\frac{\partial}{\partial x_3}\right)^{n+1}f(q)=0,\quad \forall q\in U.$$
\\
The right quaternionic vector space of polyanalytic Fueter regular functions will be denoted by $\mathcal{FR}_{n+1}(U)$.
\end{defn}
\begin{rem}
We note that two Fueter mapping theorems were proved in the polyanalytic setting, see \cite{ADS}. This will allow to introduce the so-called polyanalytic Fueter mappings $\mathcal{C}_{n+1}$ and $\tau_{n+1}$ which will be studied and investigated in the next sections.
\end{rem}
	
An important system of quaternionic polynomials that will be needed in the sequel are the so-called quaternionic Appell polynomials. This system was considered in the litterature from different points of view, for more details see \cite{CMF2017, CMF2011, DKS2019} and the references therein. Such quaternionic Appell polynomials can be defined by the following relation
\begin{equation} \label{Qk2}
Q_k(q, \bar{q})=\sum_{j=0}^k T^k_jq^{k-j}\overline{q}^j,\quad q\in\Hq, \quad k\geq 0,
\end{equation}
where 
\begin{equation}
T^k_j:=\frac{k!}{(3)_k}\frac{(2)_{k-j}(1)_j}{(k-j)!j!}=\frac{2(k-j+1)}{(k+1)(k+2)}
\end{equation}
and $(a)_n=a(a+1)...(a+n-1)$ is the Pochhammer symbol. 
\begin{rem}
Notice that the polynomials $(Q_k)_{k\geq 0}$  given by \eqref{Qk2} are Fueter regular on $\Hq$. Moreover, they form an Appell system with respect to the hypercomplex derivative $\displaystyle\frac{\overline{\mathcal{D}}}{2}:= \frac{1}{2}\left(\frac{\partial}{\partial x_0}-i\frac{\partial}{\partial x_1}-j\frac{\partial}{\partial x_2}-k\frac{\partial}{\partial x_3}\right)$. i.e., for all $k\geq 1$ we have the Appell property 
\begin{equation}
\label{a}
\frac{\overline{\mathcal{D}}}{2}Q_k(q,\overline{q})=kQ_{k-1}(q,\overline{q}).
\end{equation}
\end{rem} 
In \cite{ADS2019} the authors introduced the quaternionic polyanalytic Fock space defined for a given $I \in \mathbb{S}$ to be
$$ \mathcal{\widetilde{F}}_{Slice}^{N+1}(\mathbb{H}):= \{f \in \mathcal{SP}_{N+1}(\mathbb{H}): \int_{\mathbb{C}_I}|f_I(q)|^2 e^{-2 \pi |q|^2} d \lambda_I(q)< \infty\}, \qquad N \geq 0.$$
In \cite[Prop. 4.1]{ADS} and \cite[Prop. 4.2]{ADS} is showed that the polyanalytic Fock space is a quaternionic reproducing kernel Hilbert space which does not depend on the choice of $ I \in \mathbb{S}$. Thus, from now we will denote the quaternionic polyanalytic Fock space by $ \mathcal{\widetilde{F}}_{Slice}^{N+1}((\mathbb{H})$.
\\Now, we give the definition of the quaternionic true polyanalytic (QTP) Fock space.
\begin{defn}\label{QTPF}
A function $ f: \mathbb{H} \to \mathbb{H}$ belongs to the QTP Fock space $ \mathcal{F}^T_n (\mathbb{H})$ if and only if
\begin{itemize}
\item[i)] $ \displaystyle \int_{\mathbb{C}_I}|f_I(q)|^2 e^{-2 \pi |q|^2} \, d \lambda_I(q) < \infty.$
\item[ii)] There exists a slice regular function $H$ such that
$$ f(q)= (-1)^n\sqrt{\frac{1}{(2 \pi)^n n!}} e^{2 \pi |q|^2} \partial_s^n(e^{-2 \pi |q|^2} H(q)).$$
\end{itemize}
\end{defn}
Another important space in this context is the quaternionic Hilber space $L^{2}(\mathbb{R}; dx)=L^{2}(\mathbb{R}, \mathbb{H})$, consisting of all the square integrable $\Hq$-valued functions with respect to
\begin{align}\label{spR1}
\scal{\varphi,\psi}_{L^2(\mathbb{R};dx)} : =  \int_{\R} \overline{\psi(x)} \varphi(x) dx.
\end{align}
The real weighted Hermite functions
\begin{equation}
\label{hermi}
h_k^\nu(x) := (-1)^k e^{\frac{\nu}{2}x^2} \frac{d^k}{dx^k}\left(e^{-\nu x^2}\right)
\end{equation}

form an orthogonal basis of $L^2(\mathbb{R};dx)$, with norm given explicitly by
\begin{equation}\label{normhn}
\norm{h_k^\nu}_{L^2(\mathbb{R};dx)}^2=  2^k\nu^k k!\left(\frac{\pi}{\nu}\right)^{1/2} .
\end{equation}
Now, we recall the QTP Bargmann transform which was introduced in \cite{BEA}, and more results on this can be found on \cite{DMD2}. Let $ \varphi: \mathbb{H} \to \mathbb{R}$, then
\begin{eqnarray}
\nonumber
\nonumber
\label{pre}
(B^{n+1} \varphi)(q)&=&(-1)^n\sqrt{ \frac{1}{(2 \pi)^n n!}} e^{2 \pi |q|^2} \partial_s^n [e^{-2 \pi |q|^2} \mathcal{B}_{\mathbb{H}}\varphi(q)]\\	&=& 2^{\frac{3}{4}} (2^n n! (2 \pi)^n)^{- \frac{1}{2}} \int_{\mathbb{R}} e^{- \pi(q^2+t^2)+2 \pi \sqrt{2}qt} H_n \left( \frac{q+ \bar{q}}{\sqrt{2}}-t\right) \varphi(t) dt,
\end{eqnarray} 
where $\mathcal{B}_{\mathbb{H}}\varphi(q)=2^{\frac{3}{4}} \int_{\mathbb{R}}e^{- \pi (q^2+x^2)+\nu \sqrt{2}qx}\varphi(x) dx$ and $H_n$ are the weighted Hermite polynomials defined as
	
\begin{equation}\label{Hpi}
H_n^{2\pi}(y)=(-1)^n e^{2 \pi y^2} \frac{d^n}{dy^n} e^{-2 \pi y^2}= n! \sum_{j=0}^{\lfloor \frac{n}{2}\rfloor} \frac{(-1)^j (4 \pi y)^{n-2j}}{j!(n-2j)!},
\end{equation}

where $\lfloor.\rfloor$ denotes the integer part. 
\\In the polyanalytic setting it is possible to define a quaternionic Bargmann transform for a vector valued function $ \vec{\varphi}=(\varphi_0,..., \varphi_N)$. We will consider that a function $ \vec{\varphi}=(\varphi_0,..., \varphi_N)$ belongs to $L^2(\mathbb{R}, \mathbb{H}^{N+1})$ if and only if
$$
||\vec{\varphi}||^{2}_{L^2(\mathbb{R}, \mathbb{H}^{N+1})}:= \sum_{j=0}^{N}||\varphi_j||^{2}_{L^2(\R,\Hq)} < \infty.
$$
\begin{defn}
\label{vecba}
Let $\vec{\varphi}=(\varphi_0,...,\varphi_{N})$ be a vector-valued function in $L^{2}(\mathbb{R}, \mathbb{H}^{N+1})$. The quaternionic full-polyanalytic Bargmann (QFP) transform is defined as
\begin{equation}
\label{bar3}
\mathfrak{B} \vec{\varphi}(q)= \sum_{j=0}^{N} B^{j+1} \varphi_j(q),
\end{equation} 
where $B^{j+1} \varphi_j(q)$ is the QTP Bargmann transform, defined in \eqref{pre}. 
\end{defn}
	
A very important tool in this context are the so-called quaternionic Hermite polynomials, see \cite{EG,II,TT}
\begin{equation}\label{HermiteQ}
H^{2\pi}_{m,p}(q, \bar{q})= (2 \pi)^{p} (-1)^m e^{2 \pi |q|^2} \partial_s^m(q^pe^{-2 \pi |q|^2}), \qquad m,p \in \mathbb{N}.
\end{equation}
\begin{rem}
For any $p\geq 0$, we have 
$$H_{0,p}^{2\pi}(q,\overline{q})=(2\pi)^pq^p$$
and $$H_{1,p}^{2\pi}(q,\overline{q})=(2\pi)^{p+1}\overline{q}q^p-(2\pi)^ppq^p.$$
\end{rem}
These polynomials enjoy the following orthogonality relation, see \cite[Appendix A]{DMD2}
\begin{equation}
\label{orto}
\int_{\mathbb{C}_I} \overline{H_{m,p}^{2\pi}(q, \bar{q})} H_{m',p'}^{2\pi}(q, \bar{q}) d \lambda _I(q)= \frac{m! p! (2 \pi)^{m+p}}{2} \delta_{m, m'} \delta_{p,p'}.
\end{equation}
	
\section{The action of the polyanalytic Fueter mappings}
In this section we apply the polyanalytic Fueter maps to a slice polyanalytic function of order $n+1$. This has been done in \cite{DKS2019} for the classic Fueter map. In particular, this suggests to consider a new type of polynomials.
\subsection{The action of the map $\mathcal{C}_{n+1}$ on $ \bar{q}^k q^j$}
Let $\Omega$ be an s-domain that contains the origin and $f\in\mathcal{SP}_{n+1}(\Omega)$ be a slice polyanalytic function of order $n+1$ on $\Omega$. We note that by polyanalytic decomposition we have 
$$f(q)= \sum_{k=0}^{n}\overline{q}^kf_k(q), \quad q\in \Omega,$$ 
with $f_k\in\mathcal{SR}(\Omega)$. Then, using the series expansion theorem for slice hyperholomorphic functions we have $ f_k(q)=\sum_{j=0}^{\infty}q^j\alpha_{k,j}$ with $\lbrace{\alpha_{k,j}}\rbrace_{0 \leq k\leq n, j \geq 0}\subset\mathbb{H}$. Thus, we obtain the series expansion given by
\begin{equation}\label{series}
f(q)=\sum_{k=0}^{n}\sum_{j=0}^{\infty}\overline{q}^kq^j \alpha_{k,j}.
\end{equation}
	
We note that the polynomials $(\overline{q}^kq^j)_{j\geq 0}$ with $0\leq k  \leq n$ form the building block of our theory.
\begin{defn}{\cite[Thm. 3.12]{ADS}}
Let $\Omega$ be an axially symmetric slice domain and $f\in\mathcal{SP}_{n+1}(\Omega)$ with a polyanalytic decomposition given by $f=\displaystyle \sum_{k=0}^{n}\overline{q}^kf_k$, where $f_k\in\mathcal{SR}(\Omega)$. Then, we consider the polyanalytic Fueter map $\mathcal{C}_{n+1}:\mathcal{SP}_{n+1}(\Omega)\longrightarrow \mathcal{FR}_{n+1}(\Omega),$ which is defined by $$\displaystyle\mathcal{C}_{n+1}(f)(q):=\sum_{k=0}^nx_0^k\Delta f_k(q), \quad \text{ for any } \, q\in\Omega.$$
\end{defn}
	
\begin{prop}[Linearity of the map $\mathcal{C}_{n+1}$]
\label{lin1}
For any $f,g\in\mathcal{SP}_{n+1}(\Omega)$ and $\lambda\in\mathbb{H}$, we have $$\mathcal{C}_{n+1}(f+g\lambda)(q)=\mathcal{C}_{n+1}(f)(q)+\mathcal{C}_{n+1}(g)(q)\lambda.$$
\end{prop}
\begin{proof}
Let $f,g\in \mathcal{SP}_{n+1}(\Omega)$. We know by polyanalytic decomposition that we have $\displaystyle f(q)=\sum_{k=0}^{n}\overline{q}^kf_k(q)$ and $g(q)=\displaystyle\sum_{k=0}^{n}\overline{q}^kg_k(q)$ with $f_k, g_k\in\mathcal{SR}(\Omega)$,  for any $k=0,...,n$.  In particular, we note that  $$ (f+g\lambda)(q)=\sum_{k=0}^{n}\overline{q}^k(f_k(q)+g_k(q)\lambda).$$
Therefore \[ \begin{split}
\mathcal{C}_{n+1}(f+g\lambda)(q) & =\sum_{k=0}^{n}x_0^k\Delta(f_k+g_k\lambda )(q) \\
&= \sum_{k=0}^{n}x_0^k\Delta f_k(q)+\sum_{k=0}^{n}x_0^k\Delta g_k(q)\lambda\\\
& =\mathcal{C}_{n+1}(f)(q)+\mathcal{C}_{n+1}(g)(q)\lambda.
\\
&
\end{split}
\]
Finally, we obtain $$\mathcal{C}_{n+1}(f+g\lambda)=\mathcal{C}_{n+1}(f)+\mathcal{C}_{n+1}(g)\lambda.$$
Thus, it turns out that the map $\mathcal{C}_{n+1}$ is $\mathbb{H}$-linear.
\end{proof}
We apply the map $\mathcal{C}_{n+1}$ to the  expansion \eqref{series} and get \begin{equation}\displaystyle
\mathcal{C}_{n+1}(f):=\sum_{k=0}^{n}\sum_{j=0}^{\infty}\mathcal{C}_{n+1}(\overline{q}^kq^j)\alpha_{k,j}.
\end{equation}
So, we need to compute the action of the map $\mathcal{C}_{n+1}$ on $\overline{q}^kq^j$ with $k=0,...,n$ and $j\geq 0$. Thus, we have the following result
\begin{thm}\label{C-n action}
Let $n\geq 0$ be fixed. Then, for any $0 \leq k \leq n$ we have
\begin{equation}
\mathcal{C}_{n+1}(\bar{q}^k q^j)=
\begin{cases}
0 & \mbox{if } \quad j=0,1 \\
-2(j-1)jx_0^kQ_{j-2}(q,\overline{q}) & \mbox{if} \quad j \geq 2.
\end{cases}
\end{equation}
\end{thm}

\begin{proof}
For any $0 \leq k\leq n$, we set $f_k(q):=q^j$. Thus, we have
\[ \begin{split}
\mathcal{C}_{n+1}(\overline{q}^{k}q^j)& =\mathcal{C}_{n+1}(\overline{q}^{k}f_k)\\
&=x_0^k\Delta(f_k)(q)\\
&= x_0^k\Delta(q^j).
\end{split}
\]
We know by the proof of \cite[Prop 4.2]{DKS2019} that for $j=0,1$ we have $\mathcal{C}_{n+1}(\overline{q}^k)=\mathcal{C}_{n+1}(\overline{q}^{k}q)=0$ and for any $j\geq 2,$ we get
\[ \begin{split}
\displaystyle \mathcal{C}_{n+1}(\overline{q}^{k}q^j)& = x_0^k\Delta(q^j)\\
&=-2(j-1)jx_0^kQ_{j-2}(q,\overline{q}).\\
&
\end{split}
\]
\end{proof}
As a consequence, we note that for any $s \geq 0$ we have
	
\[ \begin{split}
\displaystyle \mathcal{C}_{n+1}(\overline{q}^kq^{s+2}) & = x_0^k\Delta(q^{s+2})\\
&=-2(s+1)(s+2)x_0^kQ_s(q,\overline{q}),\\
&
\end{split}
\]
where $Q_s(q,\overline{q})$ are defined in \eqref{Qk2}.
\newline 
This suggests to consider a new family of polynomials that we call generalized Appell polyanalytic polynomials and which are given by
\begin{equation}
\label{polmal}
\mathcal{M}_{k,s}(q,\overline{q}):=x_0^kQ_s(q,\overline{q}),\quad k=0,...,n, \quad s\geq 0.
\end{equation}
Our aim is to find a sort of Appell property for the polynomials $\mathcal{M}_{k,s}(q,\overline{q})$. In order to do this, we will first need a preliminary result
\begin{lem}\label{DjQs}
For any $s\geq j$, we have $$\overline{\mathcal{D}}^j\left(Q_s (q, \bar{q})\right)=2^j\frac{s!}{(s-j)!}Q_{s-j}(q,\bar{q}), \quad  \forall q\in \mathbb{H}.$$
\end{lem}
\begin{proof}
We prove the statement by induction. First, for $j=1$ we have
$$\overline{\mathcal{D}}\left(Q_s(q,\overline{q})\right)=2sQ_{s-1}(q,\bar{q})=2\frac{s!}{(s-1)!}Q_{s-1}(q,\overline{q}),$$
this holds true by formula \eqref{a}. Now, we suppose the result holds for $j\geq 1$ and let us prove it for $j+1$. Indeed, we have
\[ \begin{split}
\overline{\mathcal{D}}^{j+1}\left(Q_s(q,\bar{q})\right) & = \overline{\mathcal{D}} \,  \overline{\mathcal{D}}^{j}\left(Q_s(q,\bar{q})\right)\\
&=\overline{\mathcal{D}}\left(2^j\frac{s!}{(s-j)!}Q_{s-j}(q,\bar{q})\right) \\
&= 2^j\frac{s!}{(s-j)!}2(s-j)Q_{s-j-1}(q,\bar{q})
\\
&=2^j\frac{s!}{(s-j)(s-j-1)!}2(s-j)Q_{s-j-1}(q,\bar{q})
\\
&=2^{j+1}\frac{s!}{(s-j-1)!}Q_{s-j-1}(q,\bar{q}).
\end{split}
\]
\end{proof}
\begin{thm}\label{Poly-Appell}
Let $n\geq 0$ fixed, $0 \leq k \leq n$ and $s \geq k+1$. Then we have 
\begin{equation}
\overline{\mathcal{D}}^{k+1}(\mathcal{M}_{k,s}(q,\overline{q}))=\sum_{j=1}^{k+1} 2^j {k+1\choose j}\frac{k!s!}{(j-1)!(s-j)!}\mathcal{M}_{j-1,s-j}(q,\overline{q}).
\end{equation}
\end{thm}
\begin{proof}
We note that the function $x_0^k$ is real valued. Thus, we can use the Leibniz rule for the operator $\overline{\mathcal{D}}$ combined with Lemma \ref{DjQs}  and get
\[ \begin{split}
\displaystyle \overline{\mathcal{D}}^{k+1}(x_0^kQ_s(q,\bar{q})) &= \sum_{j=0}^{k+1}  {k+1\choose j}\overline{\mathcal{D}}^{k+1-j}(x_0^k)\overline{\mathcal{D}}^j Q_s(q,\bar{q})\\
&= \sum_{j=1}^{k+1}  {k+1\choose j}\frac{k!}{(k-1-k+j)!}x_{0}^{k-1-k+j}2^j\frac{s!}{(s-j)!}Q_{s-j}(q,\bar{q})\\
&=  \sum_{j=1}^{k+1} 2^j  {k+1\choose j}\frac{k!}{(j-1)!}\frac{s!}{(s-j)!}x_{0}^{j-1}Q_{s-j}(q,\bar{q})
\\
&=\sum_{j=1}^{k+1} 2^j  {k+1\choose j}\frac{k!}{(j-1)!}\frac{s!}{(s-j)!}\mathcal{M}_{j-1,s-j}(q,\bar{q}).
\end{split}
\]
\end{proof}
	
\begin{cor}
\label{appe1}
For any $n\geq 0$ fixed and $s \geq n+1$, we have
		
$$\overline{\mathcal{D}}^{n+1}(\mathcal{M}_{n,s}(q,\overline{q}))=\displaystyle\sum_{j=1}^{n+1} 2^j {n+1\choose j}\frac{n!s!}{(j-1)!(s-j)!}\mathcal{M}_{j-1,s-j}(q,\overline{q}). $$
\end{cor}
\begin{proof}
It is enough to take $k=n$ in Theorem \ref{Poly-Appell}.
\end{proof}
\begin{rem}
We note that the result proved in Corollary \ref{appe1} is an extension of the classical Appell property for the generalized Appell polyanalytic polynomials $\mathcal{M}_{k,s}(q,\bar{q})$. Indeed for $n=0$ in Corollary \ref{appe1} and by the observation that $\mathcal{M}_{0,s}(q,\bar{q})=Q_s(q,\bar{q})$ we get $$\overline{\mathcal{D}}(Q_s(q,\bar{q}))=2\frac{s!}{(s-1)!}\mathcal{M}_{0,s-1}(q,\overline{q})=2sQ_{s-1}(q,\bar{q}),\quad s\geq 1.$$
\end{rem}
\begin{rem}
The polynomials $\mathcal{M}_{k,s}$ are Fueter polyanalytic of order $n+1$. 
\end{rem}
	
\begin{rem}
\label{pm0}
We can write the action of the map $\mathcal{C}_{n+1}$ on $\overline{q}^kq^j$ in the following way when $j\geq 2$,
\begin{eqnarray}
\nonumber
\mathcal{C}_{n+1}(\overline{q}^kq^j) &=& -2(j-1)jx_0^kQ_{j-2}(q,\bar{q})\\
\label{pm}
&=& -2(j-1)j\mathcal{M}_{k,j-2}(q,\bar{q}).
\end{eqnarray}
		
\end{rem}
\begin{thm}
Let $\Omega$ be an axially symmetric slice domain. We have that $f\in\mathcal{C}_{n+1}(\mathcal{SP}_{n+1}(\Omega))$ if and only if we have the expansion 
\begin{equation}
f(q)= \sum_{k=0}^{n}\sum_{s=0}^{\infty}\mathcal{M}_{k,s}(q,\bar{q})\beta_{k,s}, \quad  \{\beta_{k,s}\}_{0 \leq k \leq n, \, s \geq 0} \subset \mathbb{H}.
\end{equation}
\end{thm}
\begin{proof}
Let $f\in\mathcal{C}_{n+1}(\mathcal{SP}_{n+1}(\Omega))$. Then, there exist a function $g\in\mathcal{SP}_{n+1}(\Omega)$ such that $$f=\mathcal{C}_{n+1}(g).$$
Thus, we know that we can write an expansion of $g$ as
$$g(q)= \sum_{k=0}^{n}\sum_{j=0}^{\infty} \overline{q}^kq^j \alpha_{k,j}, \quad \{\alpha_{k,s}\}_{0 \leq k \leq n, \, j \geq 0} \subset \mathbb{H}.$$
In particular, this implies that we have 
\begin{equation}
f(q)=\displaystyle \sum_{k=0}^{n}\sum_{j=0}^{\infty}\mathcal{C}_{n+1}(\bar{q}^kq^j)\alpha_{k,j}.
\end{equation}
Therefore, using formula \eqref{pm} we obtain 
\[ \begin{split}
\displaystyle f(q) &= -2\sum_{k=0}^n\sum_{j=2}^{\infty}(j-1)j\mathcal{M}_{k,j-2}(q,\bar{q})\alpha_{k,j}\\
&= -2\sum_{k=0}^{n}\sum_{s=0}^{\infty}(s+1)(s+2)\mathcal{M}_{k,s}(q,\bar{q})\alpha_{k,s+2}. \\
&
\end{split}
\]
Now, we set $$\displaystyle\beta_{k,s}:=-2(s+1)(s+2)\alpha_{k,s+2},$$  and  obtain  $$\displaystyle f(q)=\sum_{k=0}^{n}\sum_{s=0}^{\infty} \mathcal{M}_{k,s}(q,\bar{q})\beta_{k,s}.$$
\end{proof}
	
\begin{rem}
Now, we can observe that the map $$\mathcal{C}_{n+1}:\mathcal{SP}_{n+1}(\Omega)\longrightarrow \mathcal{FR}_{n+1}(\Omega)$$ can be introduced such that for any $f(q)=\displaystyle \sum_{k=0}^{n}\sum_{s=0}^{\infty}\overline{q}^kq^s\alpha_{k,s}$ we have $$ \mathcal{C}_{n+1}(f)(q)=\sum_{k=0}^{n}\sum_{s=0}^{\infty}\mathcal{M}_{k,s}(q,\bar{q})\beta_{k,s},$$
with
$$\displaystyle\beta_{k,s}:=-2(s+1)(s+2)\alpha_{k,s+2}.$$
\end{rem}

\subsection{The action of the polyanalytic Fueter map $\tau_{n+1}$ on $ \bar{q}^k q^j$}
We start by recalling the definition of the maps $\tau_{n+1}$, see \cite[Thm. 3.7]{ADS}.
\begin{eqnarray*}
\tau_{n+1} : \mathcal{SP}_{n+1}(\Omega) &\longrightarrow& \mathcal{FR}(\Omega)\\
f(q)= \sum_{k=0}^n \sum_{j=0}^\infty \bar{q}^k q^{j} \alpha_{kj} & \longmapsto& \tau_{n}(f)=\Delta V^n (f)(q),
\end{eqnarray*}
where $$V(f)(q)= \frac{\partial}{\partial x_0}f(q)+ \frac{\vec{q}}{| \vec{q}|^2} \sum_{\ell=1}^3 x_{\ell} \frac{\partial}{\partial x_\ell} f(q) \qquad q \in \Omega \setminus \mathbb{R}.$$
This operator is a global operator with nonconstant coefficients, and it was studied in \cite{CGS2012, CGS2013, GP}.
We will need also the following formula that puts in relation the two polyanalytic-Fueter mappings $\tau_{n+1}$ and $\mathcal{C}_{n+1}$, respectively, see .
\begin{thm}{\cite[Thm 3.13]{ADS}}
Let $f: \Omega \longrightarrow \mathbb{H}$ be a slice polyanalytic function of order $n + 1$ on an axially symmetric slice domain. Then we have
\begin{equation}
\label{R1}
\mathcal{D}^n \mathcal{C}_{n+1}(f)(q)= \frac{1}{2^n} \tau_{n+1}(f)(q),\quad \forall q\in \Omega.
\end{equation}
\end{thm}
\begin{prop}
\label{line2}
The map $\tau_{n+1}$ is $\mathbb{H}$-linear i.e.
$$ \tau_{n+1}(f+g \lambda)(q)= \tau_{n+1}(f)(q)+ \tau_{n+1}(g)(q) \lambda,$$
where $f,g \in \mathcal{SP}_{n+1}(\Omega)$ and $ \lambda \in \mathbb{H}$.
\end{prop}
\begin{proof}
From formula \eqref{R1}, the linearity of $ \mathcal{C}_{n+1}$ (see Proposition \ref{lin1}) and $ \mathcal{D}^n$ we get
\begin{eqnarray*}
\tau_{n+1}(f+g \lambda)(q) &=& 2^n \mathcal{D}^n \mathcal{C}_{n+1}(f+g \lambda) (q)\\
&=&2^n \mathcal{D}^n [\mathcal{C}_{n+1}(f)(q)+\mathcal{C}_{n+1}(g) \lambda] (q)\\
&=& 2^n [\mathcal{D}^n \mathcal{C}_{n+1}(f)(q)+\mathcal{D}^n\mathcal{C}_{n+1}(g)(q) \lambda].
\end{eqnarray*}
Using another time formula \eqref{R1} we obtain
$$ \tau_{n+1}(f+g \lambda)(q)= \tau_{n+1}(f)(q)+ \tau_{n+1}(g)(q) \lambda.$$
\end{proof}
From the series expansion \eqref{series} we have
$$ \tau_{n+1}(f)= \sum_{k=0}^n \sum_{j=0}^\infty \tau_{n+1}(\bar{q}^k q^j) \alpha_{k,j}.$$
This is the motivation for the following
\begin{thm}
\label{BL1}
For any fixed $n\geq 0$ and $0 \leq k \leq n$ we have
\begin{equation}
\tau_{n+1}(\bar{q}^k q^j)=
\begin{cases}
0 & \mbox{if } \quad j=0,1 \\
-2^{n+1} n! j(j-1) Q_{j-2}(q, \bar{q}) & \mbox{if} \quad j \geq 2.
\end{cases}
\end{equation}
\end{thm}
\begin{proof}
First of all we recall that if a function $f$ is written using the polyanalytic decomposition
$$ f(q)= \sum_{k=0}^n \bar{q}^k f_k(q),$$
then
$$ V^n f(q)=2^n n! f_{n}(q),$$
where $f_n(q)$ is the last term of the polyanalytic decomposition, see \cite{ADS}.
\\ Now, we set $f_j(q):=q^j$. By the previous property of the operator $V^n$ we have
$$ V^{n}\left(\bar{q}^k f_j(q)\right)= 2^n n! f_{j}(q)=2^n n! q^j.$$
Therefore for $j \geq 2$ we have
\begin{eqnarray*}
\tau_{n+1}(\bar{q}^k q^j) &=& \Delta V^n \left(\bar{q}^k f_j(q)\right) \\
&=& 2^n n! \Delta q^j \\
&=&-2^{n+1} n! (j-1)j Q_{j-2}(q, \bar{q}).
\end{eqnarray*}
Since $ \Delta q^0=\Delta q=0$ for the cases $j=0$ and $j=1$ we have
$$ \tau_{n+1}(\bar{q}^k)=\tau_{n+1}(\bar{q}^k q)=0.$$
\end{proof}
\begin{rem}
If we put $n=0$ in Theorem \ref{C-n action} and Theorem \ref{BL1} we get \cite[Thm. 3.2]{DKS2019}.
\end{rem}
\begin{rem}
In this case the formula for $ \tau_{n+1}(\bar{q}^k q^j)$ does not give any suggestions about new generalized Appell polyanalytic polynomials, as it happens for the computations of $ \mathcal{C}_{n+1}(\bar{q}^k q^j)$, see Theorem \ref{C-n action}. However, this behaviour is natural because $\tau_{n+1}$ maps slice polyanalytic functions of order $n+1$ in the space of Fueter regular functions, which is a well-known space. On the other side, the map $\mathcal{C}_{n+1}$ maps slice polyanalytic functions of order $n+1$ in the space of polyanalytic Fueter regular functions of order $n+1$, which is a different space with respect to the classical one and for this reason new kind of polynomials appear.
\end{rem}
\begin{thm}
\label{R5}
Let $\Omega$ be an axially symmetric slice domain. Then, a function $f$ belongs to $\tau_{n+1}\left(\mathcal{SP}_{n+1}(\Omega)\right)$ if and only if
$$ f(q)= \sum_{s=0}^\infty Q_s(q, \bar{q}) \beta_{k,s+2},$$
where $\beta_{k,s+2}:=-\sum_{k=0}^n 2^{n+1} n! (s+2)(s+1) \alpha_{k,s+2}.$
\end{thm}
\begin{proof}
Let us assume that $f\in\tau_{n+1}\left(\mathcal{SP}_{n+1}(\Omega)\right)$. Then there exists a function $g \in \mathcal{SP}_{n+1}(\Omega)$ such that
$$ f(q)= \tau_{n+1}(g)(q).$$
From the same arguments to obtain formula \eqref{series} we can write the function $g(q)$ in the following way
$$ g(q)= \sum_{k=0}^n \sum_{j=0}^\infty \bar{q}^k q^j \alpha_{k,j}.$$
Then by Theorem \ref{BL1} we get
\begin{eqnarray*}
f(q) &=& \sum_{k=0}^n \sum_{j=0}^\infty \tau_{n+1}(\bar{q}^k q^j) \alpha_{k,j}  \\
&=& -2^{n+1} n!\sum_{k=0}^n \sum_{j=2}^\infty j(j-1) Q_{j-2}(q, \bar{q}) \alpha_{k,j}\\
&=&-2^{n+1} n!\sum_{k=0}^n \sum_{s=0}^\infty (s+2)(s+1) Q_{s}(q, \bar{q}) \alpha_{k,s+2} \\
&=& \sum_{s=0}^\infty Q_{s}(q, \bar{q}) \beta_{k,s+2},
\end{eqnarray*}
where $ \beta_{k,s+2}:=-\sum_{k=0}^n 2^{n+1} n! (s+2)(s+1) \alpha_{k,s+2}.$
\end{proof}
\begin{rem}
The main difference between Theorem \ref{R5} and the results proved in \cite{ADS, ADS1} is that the coefficients of the sereis are a sum. On the other side, the only point in common is the presence of the generalized Appell polynomials $Q_s(q, \bar{q})$.
\end{rem}
The relation between the two polyanalytic Fueter mappings $ \mathcal{C}_{n+1}$ and $ \tau_{n+1}$ (see formula \eqref{R1}) gives us the following result
\begin{prop}
\label{res1}
For any fixed $n\geq 0$ and $j \geq 2$ we have
\begin{equation}
\mathcal{D}^n( \mathcal{M}_{k, j-2}(q, \bar{q}))=
\begin{cases}
n! Q_{j-2}(q, \bar{q}) & \mbox{if } \quad k=n \\
0 & \mbox{if} \quad 0 \leq k  <n.
\end{cases}
\end{equation}
\end{prop}
\begin{proof}
We start by applying the relation \eqref{R1} to $\bar{q}^k q^j$ and so we obtain
\begin{equation}
\label{R2}
\tau_{n+1}(\bar{q}^k q^j)=2^n \mathcal{D}^n \mathcal{C}_{n+1}(\bar{q}^k q^j).
\end{equation}
Now, we focus only on $k=n$. By Theorem \ref{BL1} and Theorem \ref{C-n action}, for $j \geq 2$, we have
$$ \tau_{n+1}(\bar{q}^n q^j)=-2^{n+1} n! j(j-1) Q_{j-2}(q, \bar{q}),$$
$$ \mathcal{C}_{n+1}( \bar{q}^n q^j)=-2(j-1)j x_0^n Q_{j-2}(q, \bar{q}).$$
By putting these relations in \eqref{R2} we obtain
$$ -2^{n+1} n! j(j-1) Q_{j-2}(q, \bar{q})=-2^n 2(j-1)j \mathcal{D}^n \left( x_0^n Q_{j-2}(q, \bar{q}) \right).$$
After some simplifications we get
$$ \mathcal{D}^n \left(x_0^n Q_{j-2}(q, \bar{q})\right)= n! Q_{j-2}(q, \bar{q}).$$
By the definition of generalized Appell polyanalytic polynomials and taking into account the fact that $k=n$ we get
\begin{equation}
\label{R4}
\mathcal{D}^n \left( \mathcal{M}_{k,j-2}(q, \bar{q}) \right)= n! Q_{j-2}(q, \bar{q}).
\end{equation}
Now, it remains to prove that
\begin{equation}
\label{R3}
\mathcal{D}^n \left( \mathcal{M}_{k,j-2}(q, \bar{q}) \right)=0, \qquad \forall 0 \leq k < n.
\end{equation}
We prove it by induction on $n$.
\\For $n=1$ (which means $k=0$) we have
$$ \mathcal{D}\left( \mathcal{M}_{0,j-2}(q, \bar{q}) \right)= \mathcal{D}\left( Q_{j-2}(q, \bar{q}) \right)=0.$$
We can justify the last equality from the fact that the generalized Appell polynomials are Fueter regular functions.
\\ We suppose that \eqref{R3} is true for $n$ and we want to prove it for $n+1$. Therefore we want to show that
$$ \mathcal{D}^{n+1}\left( \mathcal{M}_{k,j-2}(q, \bar{q}) \right)=0, \qquad \forall  \, \, 0 \leq k \leq n.$$
We observe that by the inductive hypothesis $ \mathcal{D}^n \left( \mathcal{M}_{k,j-2}(q, \bar{q}) \right)=0$ for all $0 \leq k <n$, thus we have
$$ \mathcal{D}^{n+1}\left(\mathcal{M}_{k,j-2}(q, \bar{q}) \right)= \mathcal{D}\mathcal{D}^{n}\left(\mathcal{M}_{k,j-2}(q, \bar{q}) \right)=0.$$
To end the proof we have only to consider the case $k=n$, this means that we have to show
$$ \mathcal{D}^{n+1}\left( \mathcal{M}_{n,j-2}(q, \bar{q}) \right)=0.$$
From \eqref{R4} and the fact that the generalized Appell polyanalytic polynomials are Fueter regular we obtain
\begin{eqnarray*}
\mathcal{D}^{n+1}\left( \mathcal{M}_{n,j-2}(q, \bar{q}) \right)&=& \mathcal{D} \mathcal{D}^n\left( \mathcal{M}_{n,j-2}(q, \bar{q}) \right)\\
&=& \mathcal{D} \left(n! Q_{j-2}(q, \bar{q}) \right) \\
&=&0.
\end{eqnarray*}
\end{proof}

\section{The true polyanalytic Fueter Bargmann transforms }
As done in \cite{ACSS}, the slice hyeprholomorphic Fock space is characterized in terms of power series, in the same spirit, here, we want to characterize the QTP Fock space denoted by $\mathcal{F}_n^T(\mathbb{H})$, see Definition \ref{QTPF}. To this end, we need to recall the quaternionic Hermite polynomials $H_{n,j}^{2\pi}(q,\bar{q})$ (see \eqref{HermiteQ}) that form an orthogonal basis of $\mathcal{F}_n^T(\mathbb{H})$, see \cite{TT}. Indeed, for any $f\in\mathcal{F}_n^T(\mathbb{H})$ we have $$f(q)=\displaystyle\sum_{j=0}^{\infty}H_{n,j}^{2\pi}(q,\bar{q})\alpha_j, \quad \lbrace \alpha_j \rbrace_{j\geq 0}\subset \mathbb{H}.$$
	
\begin{prop}\label{prop1}
A function of the form $f(q)=\displaystyle\sum_{j=0}^{\infty}H_{n,j}^{2\pi}(q,\bar{q})\alpha_j$ belongs to the space $\mathcal{F}_n^T(\mathbb{H})$ if and only if  $$\displaystyle \sum_{j=0}^{\infty}(2\pi)^j j!|\alpha_j|^2<\infty.$$
\end{prop}
\begin{proof}
Let us consider $f\in\mathcal{F}_n^T(\mathbb{H})$ this means that 
$$f(q)=\sum_{j=0}^\infty H_{n,j}^{2\pi}(q, \bar{q})\alpha_j.$$
Thus, by setting $d\mu(q):=e^{-2\pi|q|^2}dA(q)$ and using the orthogonality of the quaternionic Hermite polynomials (see \eqref{orto})  we obtain 
\[ \begin{split}
||f||_{\mathcal{F}_n^T}^2 & = \int_{\mathbb{C}_I}\left(\sum_{j=0}^\infty \overline{H_{n,j}^{2\pi}(q,\bar{q})\alpha_j}\right) \left(\sum_{s=0}^\infty H_{n,s}^{2\pi}(q,\bar{q})\alpha_s\right)d\mu(q)\\
&=\sum_{j,s=0}^{\infty}\int_{\mathbb{C}_I}\overline{\alpha_j} \overline{H_{n,j}^{2\pi}(q,\bar{q})}H_{n,s}^{2\pi}(q,\bar{q})\alpha_s d\mu(q) \\
&=\sum_{j,s=0}^{\infty} \overline{\alpha_j}\int_{\mathbb{C}_I} \overline{H_{n,j}^{2\pi}(q,\bar{q})}H_{n,s}^{2\pi}(q,\bar{q}) d\mu(q) \alpha_s
\\
&=\sum_{j,s=0}^{\infty} \overline{\alpha_j}n!j!\frac{(2\pi)^{n+j}}{2}\delta_{j,s}\overline{\alpha_s}
\\
&=\frac{n!(2\pi)^n}{2}\sum_{j=0}^{\infty}(2\pi)^j j!|\alpha_j|^2<\infty.
\end{split}
\]
\end{proof}
This means that we can write the space $\mathcal{F}_n^T(\mathbb{H})$ in the following way $$\displaystyle \mathcal{F}^T_n(\mathbb{H})=\left\{f= \sum_{j=0}^\infty H_{n,j}^{2\pi}\alpha_j,\quad \sum_{j=0}^{\infty}(2\pi)^j j!|\alpha_j|^2<\infty\right \}.$$
\subsection{The map $ \mathcal{C}_{n+1}$ applied to the QTP Bargmann transform}
We start by studying the action of the map $\mathcal{C}_{n+1}$ on the quaternionic Hermite polynomials. Indeed, we first prove the following result which will be very important for the sequel
\begin{prop} \label{prop2}
For any fixed $n, j\geq 0$, it holds that 
$$ \mathcal{C}_{n+1}\left(H_{n,j}^{2\pi}(q,\bar{q})\right)=-2(2\pi)^j n!j!\sum_{s=0}^{n}\frac{(-1)^s(2\pi)^{n-s}}{s!(n-s)!(j-s-2)!}\mathcal{M}_{n-s,j-s-2}(q,\bar{q}),  \quad  \hbox{if} \quad j\geq n+2,$$
otherwise
$$ \mathcal{C}_{n+1}\left(H_{n,j}^{2\pi}(q,\bar{q})\right)=0 \quad  \hbox{if} \quad j< n+2.$$
\end{prop}
\begin{proof}
We note that it is possible to write the quaternionic Hermite polynomials as follows
$$H_{n,j}^{2\pi}(q,\bar{q})=(2\pi)^jn!\sum_{s=0}^n(-1)^s\frac{j!(2\pi)^{n-s}}{s!(n-s)!(j-s)!}\bar{q}^{n-s}q^{j-s}.$$
Now, we observe that applying the polyanalytic Fueter mapping $\mathcal{C}_{n+1}$ and using Remark \ref{pm0} we get
$$\mathcal{C}_{n+1}(\bar{q}^{n-s}q^{j-s})=-2(j-s-1)(j-s)\mathcal{M}_{n-s,j-s-2}(q,\bar{q}).$$
By linearity of the map $\mathcal{C}_{n+1}$ we obtain 
\[ \begin{split}
\mathcal{C}_{n+1}\left(H_{n,j}^{2\pi}(q,\bar{q})\right)& =(2\pi)^jn!\sum_{s=0}^{n}\frac{(-1)^sj!(2\pi)^{n-s}}{s!(n-s)!(j-s)!}\mathcal{C}_{n+1}(\bar{q}^{n-s}q^{j-s})\\
&=-2(2\pi)^jn!\sum_{s=0}^{n}\frac{(-1)^sj!(j-s-1)(j-s)(2\pi)^{n-s}}{s!(n-s)!(j-s)!}\mathcal{M}_{n-s,j-s-2}(q,\bar{q}) \\
&=-2(2\pi)^jn!\sum_{s=0}^{n}\frac{(-1)^sj!(j-s-1)(j-s)(2\pi)^{n-s}}{s!(n-s)!(j-s)(j-s-1)(j-s-2)!}\mathcal{M}_{n-s,j-s-2}(q,\bar{q}) 
\\
&=-2(2\pi)^jn!\sum_{s=0}^{n}\frac{(-1)^sj!(2\pi)^{n-s}}{s!(n-s)!(j-s-2)!}\mathcal{M}_{n-s,j-s-2}(q,\bar{q}). 
\\
\end{split}
\]
\end{proof}
We denote by $\mathcal{A}_{n+1}(\mathbb{H})$ the range of the polyanalytic Fueter mapping $\mathcal{C}_{n+1}$ on the QTP Fock space. Indeed, we have
\begin{equation}
\label{ff}
\mathcal{A}_{n+1}(\mathbb{H}):=\lbrace{\mathcal{C}_{n+1}(f), \textbf{  } f\in\mathcal{F}_n^T(\mathbb{H})}\rbrace.
\end{equation}
We prove the following characterization of the previous space 
\begin{thm}\label{Space5}
Let $n\geq 0$ fixed, it holds that
$$\mathcal{A}_{n+1}(\mathbb{H})=\left\{\sum_{h=0}^{\infty}\sum_{s=0}^n\mathcal{M}_{n-s,h+n-s}(q,\bar{q})\beta_{h,s}, \, \, \sum_{s=0}^n\frac{[(h+n-s)]^2(s!)^2[(n-s)!]^2}{(2\pi)^{h+n+2-2s}(h+n+2)!}|\beta_{h,s}|^2 < \infty\right \}, $$
where $ \{\beta_{h,s}\}_{h \geq 0, \, 0 \leq s \leq n}\subset \mathbb{H}$.
\end{thm}
\begin{proof}
Let $g\in\mathcal{A}_{n+1}(\mathbb{H}),$ then then there exists a function $f\in\mathcal{F}_{n}^{T}(\mathbb{H})$ such that we have \begin{equation}\label{4}
g(q)=\mathcal{C}_{n+1}(f)(q).
\end{equation}
We note that by using Proposition \ref{prop1} we have that 
$$f(q)= \sum_{j=0}^{\infty}H_{n,j}(q,\bar{q})\alpha_j, \quad \{\alpha_j\}_{j\geq 0}\subset\mathbb{H}$$
and $$\displaystyle ||f||^{2}_{\mathcal{F}_n^T}=\sum_{j=0}^{\infty}(2\pi)^j j!|\alpha_j|^2<\infty.$$
Then, by Proposition \ref{prop2} we have that 
\[ \begin{split}
g(q) & =-2n!\sum_{s=0}^{n}\sum_{j=n+2}^{\infty}\frac{(2\pi)^j(-1)^sj!(2\pi)^{n-s}}{s!(n-s)!(j-s-2)!}\mathcal{M}_{n-s,j-s-2}(q,\bar{q})\alpha_j \\
&=-2n!\sum_{s=0}^{n}\sum_{h=0}^{\infty}\frac{(2\pi)^{h+n+2}(-1)^s(h+n+2)!(2\pi)^{n-s}}{s!(n-s)!(h+n-s)!}\mathcal{M}_{n-s,h+n-s}(q,\bar{q})\alpha_{h+n+2}  \\
&=\sum_{h=0}^{\infty}\sum_{s=0}^n\mathcal{M}_{n-s,h+n-s}(q,\bar{q})\beta_{h,s},
\\
&
\end{split}
\]
where we have set $$\beta_{h,s}:=-\frac{2n!(2\pi)^{h+n+2}(-1)^s(h+n+2)!(2\pi)^{n-s}}{s!(n-s)!(h+n-s)!}\alpha_{h+n+2}.$$
By developing the calculations we get
\begin{eqnarray*}
\sum_{h=0}^\infty \sum_{s=0}^n\frac{[(h+n-s)]^2(s!)^2[(n-s)!]^2}{(2\pi)^{h+n+2-2s}(h+n+2)!}|\beta_{h,s}|^2 \! \! \! \! \! &=&  \! \! \! \! \! 4(n!)^2(2 \pi)^{2n}\sum_{s=0}^n\sum_{h=0}^{\infty}(2\pi)^{h+n+2}(h+n+2)!|\alpha_{h+n+2}|^2\\
\! \! \! \! \!&=&  \! \! \! \! \! 4 n! (n+1)!(2 \pi)^{2n}\sum_{h=0}^{\infty}(2\pi)^{h+n+2}(h+n+2)!|\alpha_{h+n+2}|^2\\
\! \! \! \! \!&=&  \! \!  \! \! \! \|f\|^{2}_{\mathcal{F}^{T}_{n}}<\infty.
\end{eqnarray*}

Therefore
$$\sum_{h=0}^\infty \sum_{s=0}^n\frac{[(h+n-s)]^2(s!)^2[(n-s)!]^2}{(2\pi)^{h+n+2-2s}(h+n+2)!}|\beta_{h,s}|^2    < \infty.$$
On the other hand, if we consider a function of the following form $$\displaystyle h(q)=\sum_{j=n+2}^{\infty}H_{n,j}^{2\pi}(q,\bar{q})\gamma_j, $$
where we have set 
$$\gamma_j=-\frac{s!(n-s)!(j-s-2)!}{2n!(2\pi)^{j}(-1)^s(j)!(2\pi)^{n-s}}\beta_{j-n-2,s}, \qquad j \geq n+2.$$
Then, we get $g=\mathcal{C}_{n+1}(h)$ since for $0 \leq s \leq n$, we have 
$$ \sum_{s=0}^{n}\frac{(-1)^s (h+s+2)!(2\pi)^{h+n+2} (2\pi)^{n-s}}{s!(n-s)!(h+n-s)!}\mathcal{M}_{n-s,h+n-s}(q,\bar{q})=-\frac{1}{2n!}\mathcal{C}_{n+1}(H_{n,h+n+2}^{2\pi}(q,\bar{q})).$$
Furthermore, we note that for any fixed  index $0 \leq s \leq n$ we have 
\[ \begin{split}
||h||^{2}_{\mathcal{F}_n^T} & =\sum_{j=n+2}^{\infty}(2\pi)^j j!|\gamma_{j}|^2 \\
&= \frac{1}{4(n!)^2}\sum_{\ell=0}^\infty \frac{(2\pi)^{\ell+n+2}(\ell+n+2)!(s!)^2[(n-s)!]^2(\ell+n-s !)^2}{(2\pi)^{2\ell+2n+4}(2\pi)^{2(n-s)}[(\ell+n+2)!]^2}|\beta_{\ell,s}|^2.  \\
& =\frac{(s!)^2[(n-s)!]^2}{4(n!)^2(2\pi)^{n+2}(2\pi)^{2(n-s)}}\sum_{\ell=0}^{\infty}\frac{[(\ell+n-s !)]^2}{(2\pi)^\ell (\ell+s+2)!}| \beta_{\ell,s}|^2<\infty.
\\
&
\end{split}
\]
Hence, we have $g=\mathcal{C}_{n+1}(h)$ with $h\in\mathcal{F}_{n}^{T}(\mathbb{H})$, this shows that $h\in\mathcal{A}_{n+1}(\mathbb{H})$.
\end{proof}
	
\begin{defn}
\label{inner}
Let $f,g\in\mathcal{A}_{n+1}(\mathbb{H})$ be such that $\displaystyle f=\sum_{h=0}^{\infty}\sum_{s=0}^n\mathcal{M}_{n-s,h+n-s}(q,\bar{q})\alpha_{h,s}$ and $\displaystyle g=\sum_{h=0}^{\infty}\sum_{s=0}^n\mathcal{M}_{n-s,h+n-s}(q,\bar{q})\beta_{h,s}$. Then, we define their inner product as
$$ \scal{f,g}_{\mathcal{A}_{n+1(\mathbb{H})}}:=\sum_{h=0}^\infty \sum_{s=0}^n\frac{[(h+n-s)]^2(s!)^2[(n-s)!]^2}{(2\pi)^{h+n+2-2s}(h+n+2)!}\overline{\beta_{h,s}}\alpha_{h,s}.$$
The associated norm is given by 
		
$$\displaystyle ||f||^2=\scal{f,f}_{\mathcal{A}_{n+1}(\mathbb{H})}=\sum_{h=0}^\infty \sum_{s=0}^n\frac{[(h+n-s)]^2(s!)^2[(n-s)!]^2}{(2\pi)^{h+n+2-2s}(h+n+2)!}|\alpha_{h,s}|^2.$$
\end{defn}
\begin{rem}
If we take $n=0$ in the space characterized by Theorem \ref{Space5}
we obtain the space $\mathcal{A}(\mathbb{H})$ studied in \cite{DKS2019}. Indeed, if $n=0$ we obtain
$$\mathcal{A}_1(\mathbb{H})=\left\lbrace \sum_{h=0}^{\infty}\mathcal{M}_{0,h}(q,\bar{q})\beta_{h,0}, \quad (\beta_{h,0})_{h\geq 0}\subset \mathbb{H}, \quad \sum_{h=0}^{\infty}\frac{(h!)^2}{(2\pi)^h(h+2)!}|\beta_{h,0}|^2<\infty \right\rbrace.$$
This observation holds since $\mathcal{M}_{0,h}(q, \bar{q})=Q_h(q, \bar{q})$ and $$ \sum_{h=0}^{\infty}\frac{h! h!}{(2\pi)^h(h+2)(h+1) h!}=\sum_{h=0}^{\infty}\frac{h!}{(2\pi)^h(h+2)(h+1)}.$$
Thus, we get that 
$$\mathcal{A}_1(\mathbb{H})=\mathcal{A}(\mathbb{H}).$$
\end{rem}
In the next result, we write the kernel of the QTP Bargmann transform as a generating function of the quaternionic Hermite polynomials.
\begin{prop}\label{kernelpoly}
Let $n\geq 0$ fixed. Then, for any $q\in\mathbb{H}$ and $x\in\mathbb{R}$, we have 
$$ \sum_{k=0}^{\infty}\frac{h_{k}^{2\pi}(x)}{||h_{k}^{2\pi}||_{L^2_{\mathbb{H}}(\mathbb{R})}}\frac{H_{n,k}^{2\pi}(q,\bar{q})}{||H^{2\pi}_{n,k}||_{\mathcal{F}_n^T(\mathbb{H})}}=2^{\frac{3}{4}}(2^nn!(2\pi)^n)^{-\frac{1}{2}}e^{-\pi(q^2+x^2)+2\pi\sqrt{2}qx}H_{n}^{2\pi}\left(\frac{q+\bar{q}}{\sqrt{2}}-x \right),$$
where $H_{n}^{2 \pi}$ are the Hermite polynomials defined in formula \eqref{Hpi}.
\end{prop}
\begin{proof}
We know that the QTP Bargmann transform maps $L^2(\mathbb{R},\mathbb{H})$ onto $\mathcal{F}_n^T(\mathbb{H})$. Moreover, it is well-known that the normalized Hermite function form an orthonormal basis of $L^2(\mathbb{R},\mathbb{H})$. 
We note that the normalized quaternionic Hermite polynomials $ \frac{H_{n,k}^{2\pi}(q, \bar{q})}{\|H_{n,k}^{2\pi}\|_{\mathcal{F}_n^T(\mathbb{H})}}$ form an orthonormal basis of the QTP Fock space $\mathcal{F}_n^T(\mathbb{H})$.
Now, using a classical approach used in \cite{Gaz} we know that 
\begin{equation}\label{K6}
K(q,x):= \sum_{k=0}^{\infty}\frac{h_{k}^{2\pi}(x)}{||h_{k}^{2\pi}||_{L^2_{\mathbb{H}}(\mathbb{R})}}\frac{H_{n,k}^{2\pi}(q,\bar{q})}{||H_{n,k}^{2\pi}||_{\mathcal{F}_n^T(\mathbb{H})}}
\end{equation}
is the kernel function that leads to the following integral transform
$$\displaystyle \int_{\mathbb{R}} K(q,x)f(x)dx.$$
On the other hand, we know by \cite{DMD2} that 
\begin{equation}\label{K7}
K(q,x)=2^{\frac{3}{4}}(2^nn!(2\pi)^n)^{-\frac{1}{2}}e^{-\pi(q^2+x^2)+2\pi\sqrt{2}qx}H_n\left(\frac{q+\bar{q}}{\sqrt{2}}-x \right).
\end{equation}
Putting equal the expressions \eqref{K6} and \eqref{K7} we get the thesis.
\end{proof}
	
\begin{rem}
We note that if we put $n=0$ in Proposition \ref{kernelpoly} we obtain the same result in \cite[Prop 4.1]{DG}.
\end{rem}
Now, we give the following
\begin{defn}[$ \mathcal{C}$-polyanalytic Fueter Bargmann transform]
\label{CB}
We define $ B_{\mathcal{C}}^{n+1}:L^{2}_{\mathbb{H}}(\mathbb{R}) \to \mathcal{A}_{n+1}(\mathbb{H})$ as
\begin{equation}
\label{bc}
B^{n+1}_{\mathcal{C}}:= \mathcal{C}_{n+1}\circ B^{n+1},
\end{equation}
where $B^{n+1}$ is the QTP Bargmann transform, see \eqref{pre}.
\end{defn}
	
\begin{cor}
The $ \mathcal{C}$-polyanalytic Fueter Bargmann transform $B^{n+1}_{\mathcal{C}}$ is $ \mathbb{H}$-linear.
\end{cor}
\begin{proof}
It is a direct consequence of the linearity of the map $ \mathcal{C}_{n+1}$ (see Proposition \ref{lin1}).
\end{proof}
Now, we have all what we need to write an expression of the $ \mathcal{C}$-polyanalytic Fueter Bargmann transform 
	
\begin{thm}
\label{CB1}
The $ \mathcal{C}$-polyanalytic Fueter Bargmann transform can be realized using the following commutative diagram 
$$B_{\mathcal{C}}^{n+1}: \xymatrix{L^2_\Hq(\R) \ar[r] \ar[d]_{B^{n+1}} & \mathcal{A}_{n+1}(\Hq)  \\ \mathcal{F}^{n}_{T}(\Hq) \ar[r]_{Id} & \mathcal{SP}_{n+1}(\Hq) \ar[u]_{\mathcal{C}_{n+1}}
}$$
More precisely, for any $\varphi\in L^2_\mathbb{H}(\mathbb{R})$ and $q\in\mathbb{H}$ we have $$\displaystyle B^{n+1}_{\mathcal{C}}[\varphi](q)=\int_{\mathbb{R}}\Phi(q,x)\varphi(x)dx, $$
where 
\begin{equation}\label{Phi}
\Phi(q,x):=-\sqrt{n!}(2\pi)^{\frac{n}{2}}2^{\frac{3}{4}}\sum_{l=0}^{\infty}\sum_{j=0}^{n}\frac{(-1)^j h_{l+n+2}^{2\pi}(x)\mathcal{M}_{n-j,l+n-j}(q,\bar{q})}{2^{\frac{l+n}{2}}(2\pi)^j j!(l+n-j)!(n-j)!}.
\end{equation}
and $\mathcal{M}_{n-j,l+n-j}$ are polynomials defined in \eqref{polmal}.
\end{thm}
\begin{proof}
According to Proposition \ref{kernelpoly} we know that we can write the kernel of the QTP Bargmann transform as $$K(q,x)=\displaystyle \sum_{k=0}^{\infty}\frac{h_{k}^{2\pi}(x)}{||h_{k}^{2\pi}||_{L^2_{\mathbb{H}}(\mathbb{R})}} \frac{H_{n,k}^{2\pi}(q,\bar{q})}{\|H_{n,k}(q, \bar{q})\|_{\mathcal{F}_n^T(\mathbb{H})}}.$$
We develop the computations using the explicit expressions of $||h_{k}^{2\pi}||_{L^2_{\mathbb{H}}(\mathbb{R})}$ and $||H_{n,k}||_{\mathcal{F}_n^T(\mathbb{H})}$ together with Proposition \ref{prop2} we have 
\begin{eqnarray*}
\mathcal{C}_{n+1}\left(K(q,x)\right) & =& \sum_{k=n+2}^{\infty}\frac{h_{k}^{2\pi}(x)  \, \mathcal{C}_{n+1}\left(H_{n,k}^{2\pi}(q,\bar{q})\right)}{\|h_{k}^{2\pi}\|_{L^2_{\mathbb{H}}(\mathbb{R})}\|H_{n,k}\|_{\mathcal{F}_n^T(\mathbb{H})}}\\
%&=& -\sum_{j=0}^{n}\sum_{k=n+2}^{\infty} \frac{h_k^{2 \pi}(x)2 (2 \pi)^k n! k! (2 \pi)^{n}(-1)^j (2 \pi)^{-j}}{2^{\frac{k}{2}} (2 \pi)^{\frac{k}{2}} \sqrt{k!} 2^{- \frac{1}{4}} 2^{- \frac{1}{2}} \sqrt{k!} \sqrt{n!} (2 \pi)^{\frac{k+n}{2}}j!(n-j)!(k-j-2)!}\mathcal{M}_{n-j, k-j-2}(q, \bar{q})\\
&=&-2\sqrt{n!}(2\pi)^{\frac{n}{2}}2^{\frac{3}{4}}\sum_{j=0}^{n}\sum_{k=n+2}^{\infty}\frac{(-1)^jh_{k}^{2\pi}(x)}{2^{\frac{k}{2}}(2\pi)^j j!(n-j)!(k-j-2)!}\mathcal{M}_{n-j,k-j-2}(q,\bar{q}) \\
&=&-\sqrt{n!}(2\pi)^{\frac{n}{2}}2^{\frac{3}{4}}\sum_{l=0}^{\infty}\sum_{j=0}^{n}\frac{(-1)^jh_{l+n+2}^{2\pi}(x)\mathcal{M}_{n-j,l+n-j}(q,\bar{q})}{2^{\frac{l+n}{2}}(2\pi)^j j!(l+n-j)!(n-j)!}.
\end{eqnarray*}
		
This ends the proof.
\end{proof}
\begin{rem}
If we put $n=0$ in formula \eqref{Phi} we obtain $$ \Phi(q,x)=-2^{\frac{3}{4}}\sum_{k=0}^{\infty}\frac{Q_k(q,\bar{q})h_{k+2}^{2\pi}(x)}{k!2^{\frac{k}{2}}}, \quad \forall (q,x)\in\mathbb{H}\times \mathbb{R}.$$
This is the same formula obtained in \cite[Thm 4.14]{DKS2019} up to a constant.
\end{rem}
Now we study further properties of the transform $B_{\mathcal{C}}^{n+1}$. Firstly, we need the following result
	
\begin{prop}\label{action}
Let $n\geq 0,$ we set $ \psi_{k}^{2\pi}(x)=\frac{h_{k}^{2\pi}(x)}{||h^{2\pi}_{k}||_{L^2_\mathbb{H}(\mathbb{R})}}$. Then, we have $$ B^{n+1}(\psi_{k}^{2\pi})(q)=\sqrt{\frac{2}{(2\pi)^{n+k}n!k!}}H_{n,k}^{2\pi}(q,\bar{q}).$$
\end{prop}
\begin{proof}
We note that in \cite[Lemma 4.4]{DG} we have 
$$ B_{\mathbb{H}}(\psi_{k}^{2\pi})(q)=\sqrt{2}\frac{(2\pi)^{\frac{k}{2}}}{\sqrt{k!}}q^k.$$
Thus, using the definition of the QTP Bargmann transform we get
		
\[ \begin{split}
\displaystyle B^{n+1}(\psi_{k}^{2\pi})(q) & =(-1)^n\sqrt{\frac{1}{(2\pi)^n n!}}e^{2\pi |q|^2}\partial_{s}^{n}\left(e^{-2\pi|q|^2}B_{\mathbb{H}}(\psi_{k}^{2\pi})(q)\right)\\
&=(-1)^n\sqrt{\frac{1}{(2\pi)^n n!}}e^{2\pi |q|^2}\partial_{s}^{n}\left(e^{-2\pi|q|^2} \sqrt{2}\frac{(2\pi)^{\frac{k}{2}}}{\sqrt{k!}}q^k \right) \\
&=\sqrt{\frac{1}{(2\pi)^n n!}}\frac{\sqrt{2}}{\sqrt{(2\pi)^k k!}}  (-1)^n(2\pi)^k e^{2\pi |q|^2} \partial_{s}^{n}\left(e^{-2\pi|q|^2} q^k \right) \\
&= \sqrt{\frac{1}{(2\pi)^n n!}}\frac{\sqrt{2}}{\sqrt{(2\pi)^k k!}}  H_{n,k}^{2\pi}(q,\bar{q})\\
&=\sqrt{\frac{2}{(2\pi)^{k+n}n!k!}}H_{n,k}^{2\pi}(q,\bar{q}).
\end{split}
\] 
This ends the proof.
\end{proof}
\begin{rem}
If we put $n=0$ in Proposition \ref{action} we obtain the same result of \cite[Lemma 4,4]{DG}.
\end{rem}
\begin{prop}\label{prop5}
For any fixed $n\geq 0$ we have 
$$
B^{n+1}_{\mathcal{C}}(\psi_{k}^{2\pi})(q)=
\begin{cases}
-2^{\frac{3}{2}}\sqrt{(2\pi)^{k+n}n!k!}\sum_{j=0}^{n}\frac{(-1)^j\mathcal{M}_{n-j,k-j-2}(q,\bar{q})}{(2\pi)^j j!(n-j)!(k-j-2)!} & \mbox{if } \quad k\geq n+2 \\
0, & \mbox{if} \quad k< n+2.
\end{cases}
$$
\end{prop}
\begin{proof}
We apply Proposition \ref{prop2} and Proposition \ref{action} to get 
\[ \begin{split}
\displaystyle B_{\mathcal{C}}^{n+1}(\psi_{k}^{2\pi})(q,\bar{q}) & =\sqrt{\frac{2}{(2\pi)^{k+n}n!k!}}\mathcal{C}_{n+1}(H_{n,k}^{2\pi}(q,\bar{q}))\\
&=-\sqrt{\frac{2}{(2\pi)^{k+n}n!k!}}(2 (2\pi)^kn!k!) \sum_{j=0}^{n}\frac{(-1)^j(2\pi)^{n-j}}{j!(n-j)!(k-j-2)!}\mathcal{M}_{n-j,k-j-2}(q,\bar{q}) \\
&=-2^{\frac{3}{2}}\sqrt{\frac{1}{(2\pi)^{k+n}n!k!}} (2\pi)^{k+n}n!k!\sum_{j=0}^{n}\frac{(-1)^j(2\pi)^{-j}}{j!(n-j)!(k-j-2)!}\mathcal{M}_{n-j,k-j-2}(q,\bar{q})
\\
&=-2^{\frac{3}{2}}\sqrt{(2\pi)^{k+n}n!k!} \sum_{j=0}^{n}\frac{(-1)^j}{j!(n-j)!(k-j-2)!(2\pi)^{j}}\mathcal{M}_{n-j,k-j-2}(q,\bar{q}). 
\\
\end{split}
\]
\end{proof}
	
\begin{rem}
If we put $n=0$ in Proposition \ref{prop5} we get 
$$
B^1_{\mathcal{C}}(\psi_{k}^{2\pi})(q)=
\begin{cases}
-2^{\frac{3}{2}}\frac{\sqrt{(2\pi)^k k!}}{(k-2)!}Q_{k-2}(q,\bar{q}) & \mbox{if } \quad k\geq 2 \\
0 & \mbox{if} \quad k=0,1
\end{cases}
$$
which is exactly the same result of \cite[Prop. 4.18]{DKS2019}, up to a constant.
\end{rem}
\begin{defn}
\label{hil}
Let us consider the following subspace of $L^2_{\mathbb{H}}(\mathbb{R})$
$$ \mathcal{H}^n:= \bigoplus_{k=n+2}^\infty \{ \psi_k^{2 \pi} \alpha, \, \, \alpha \in \mathbb{H} \},$$
where $ \psi_k^{2 \pi}$ are the normalized Hermite functions.
\end{defn}

\begin{rem}
Since $ \psi_k^{2 \pi}$ is an orthogonal basis for $L^2_{\mathbb{H}}(\mathbb{R})$ it is clear that
$$ L^2_{\mathbb{H}}( \mathbb{R})=\bigoplus_{k=0}^\infty  \, \, \{ \psi_k^{2 \pi} \alpha, \, \, \alpha \in \mathbb{H} \}.$$
\end{rem}
\begin{prop}
\label{split}
Let $ \varphi \in \mathcal{H}^n$. Then for any fixed $n \geq 0$ we have
$$B^{n+1}_{\mathcal{C}} \varphi= \sum_{h=0}^\infty \sum_{j=0}^n \mathcal{M}_{n-j,h+n-j}(q, \bar{q}) \alpha_{h,j} $$
where
$$\alpha_{h,j}:=-2^\frac{3}{2} \sqrt{(2 \pi)^{h+2n+2} n! (h+n+2)!} \frac{(-1)^j \alpha_{h+n+2}}{(2 \pi)^j j! (n-j)! (h+n-j)!}, \qquad \{\alpha_h\}_{h \geq 0} \subset \mathbb{H}.$$
\end{prop}
\begin{proof}
Let $ \varphi \in \mathcal{H}^n$. By Definition \ref{hil} we can write
$$ \varphi= \sum_{k=n+2}^\infty \psi_k^{2 \pi} \alpha_k, \qquad \{\alpha_h\}_{k \geq 0} \subset \mathbb{H}.$$
By Proposition \ref{prop5} we have
\begin{eqnarray*}
B^{n+1}_{\mathcal{C}} \varphi (q) \! \! \! \! \! \!&=&  \sum_{k=n+2}^{\infty}  B^{n+1}_{\mathcal{C}} (\psi_k^{2 \pi} )(q) \alpha_k \\
\nonumber
&=& -2^{\frac{3}{2}} \sum_{j=0}^n \sum_{k=n+2}^{\infty} \sqrt{(2 \pi)^{k+n} n!  k!} \frac{(-1)^j \mathcal{M}_{n-j,k-j-2}(q, \bar{q})}{(2 \pi)^j j! (n-j)! (k-j-2)!} \alpha_k\\ \nonumber
&=& -2^{\frac{3}{2}} \sum_{h=0}^\infty \sum_{j=0}^n \sqrt{(2 \pi)^{h+2n+2} n! (h+n+2)!} \frac{(-1)^j \mathcal{M}_{n-j,h+n-j}(q, \bar{q})}{(2 \pi)^j j! (n-j)! (h+n-j)!} \alpha_{h+n+2}\\ \nonumber
&:=&\sum_{h=0}^\infty \sum_{j=0}^n \mathcal{M}_{n-j,h+n-j}(q, \bar{q}) \alpha_{h,j}, \nonumber
\end{eqnarray*}
\end{proof}
\begin{prop}
\label{ris3}
For any fixed $n\geq 0$, the $\mathcal{C}$-polyanalytic Fueter Bargmann transform is a quaternionic right linear bounded surjective operator from $ L^{2}_{\mathbb{H}}(\mathbb{R})$ to $ \mathcal{A}_{n+1}(\mathbb{H})$. This means that
$$ \| B^{n+1}_{\mathcal{C} }\varphi \|_{\mathcal{A}_{n+1}(\mathbb{H})} \leq 2 \sqrt{2} \sqrt{(n+1)!} (2 \pi)^\frac{n}{2}  \| \varphi \|_{L^{2}_{\mathbb{H}}(\mathbb{R})}.$$
\end{prop}
\begin{proof}
By Proposition \eqref{split} we have that
$$ B^{n+1}_{\mathcal{C}} \varphi= \sum_{h=0}^\infty \sum_{j=0}^n \mathcal{M}_{n-j,h+n-j}(q, \bar{q}) \alpha_{h,j}',$$
with
 $$ \alpha_{h,j}':=-2^\frac{3}{2} \sqrt{(2 \pi)^{h+2n+2} n! (h+n+2)!} \frac{(-1)^j \alpha_{h+n+2}}{(2 \pi)^j j! (n-j)! (h+n-j)!}.$$
From the definition of norm of the space $ \mathcal{A}_{n+1}(\mathbb{H})$ (see Definition \ref{inner}) we get
\begin{eqnarray*}
\| B^{n+1}_{\mathcal{C}} \varphi \|_{\mathcal{A}_{n+1}( \mathbb{H})}^2 &=& \sum_{h=0}^\infty \sum_{j=0}^n\frac{[(h+n-j)]^2(j!)^2((n-j)!)^2}{(2\pi)^{h+n+2-2j}(h+n+2)!}    |\alpha_{h,j}'|^2  \\
&=&  \sum_{h=0}^{\infty}\sum_{j=0}^n 8 \frac{[(h+n-j)]^2(j!)^2((n-j)!)^2}{(2\pi)^{h+n+2-2j}(h+n+2)!} \cdot\\ 
&& \cdot  \frac{(2 \pi)^{h+2n+2} n! (h+n+2)!}{(2 \pi)^{2j} [j!(n-j)!]^2 [(h+n-j)!]^2} | \alpha_{h+n+2}|^2 \\
&=& 8 n! (2 \pi)^n \sum_{j=0}^n\sum_{h=0}^{\infty} | \alpha_{h+n+2}|^2\\
&=&  8 n! (2 \pi)^n \sum_{j=0}^n \sum_{\ell=n+2}^{\infty} | \alpha_{\ell}|^2\\
&\leq&  8 n! (2 \pi)^n \sum_{j=0}^n \sum_{\ell=0}^{\infty} | \alpha_{\ell}|^2 \\
&=& 8 n!(n+1) (2 \pi)^n \| \varphi \|_{L^{2}_{\mathbb{H}}(\mathbb{R})}^2\\
&=& 8 (n+1)! (2 \pi)^n \| \varphi \|_{L^{2}_{\mathbb{H}}(\mathbb{R})}^2.
\end{eqnarray*}
\end{proof}

\begin{prop}[Isometry of $B^{n+1}_{\mathcal{C}}$]
\label{ris2}
Let $n\geq 0$ fixed and assume that $ \varphi, \psi \in \mathcal{H}^n$. Then we have
\begin{equation}
\label{iso}
\langle B^{n+1}_{\mathcal{C}} \varphi, B^{n+1}_{\mathcal{C}} \psi \rangle_{\mathcal{A}_{n+1}(\mathbb{H})}=  8 (n+1)! (2 \pi)^n \langle \varphi, \psi \rangle_{ \mathcal{H}^n} .
\end{equation}
In particular
$$ \| B^{n+1}_{\mathcal{C}} \varphi \|_{\mathcal{A}_{n+1}(\mathbb{H})}= 2 \sqrt{2}\sqrt{(n+1)!} (2 \pi)^{\frac{n}{2}} \| \varphi \|_{\mathcal{H}^n}.$$
\end{prop}
\begin{proof}
Let us consider
$$ \varphi=\sum_{k=n+2}^\infty \psi_k^{2 \pi}\alpha_k, \qquad \psi=\sum_{k=n+2}^\infty \psi_k^{2 \pi} \alpha_k',$$
with $ \{\alpha_k\}_{k\geq 0} \subset \mathbb{H}$ and $ \{ \alpha_k'\}_{k\geq 0}\subset \mathbb{H}$
By Proposition \ref{split} we have
$$
B^{n+1}_{\mathcal{C}} \varphi=\sum_{h=0}^\infty \sum_{j=0}^n \mathcal{M}_{n-j,h+n-j}(q, \bar{q}) \alpha_{h,j}', \qquad B^{n+1}_{\mathcal{C}} \psi=\sum_{h=0}^\infty \sum_{j=0}^n \mathcal{M}_{n-j,h+n-j}(q, \bar{q}) \beta_{h,j}',
$$
where
$$ \alpha_{h,j}':=-2^\frac{3}{2} \sqrt{(2 \pi)^{h+2n+2} n! (h+n+2)!} \frac{(-1)^j \alpha_{h+n+2}}{(2 \pi)^j j! (n-j)! (h+n-j)!}, $$
$$ \beta_{h,j}':=-2^\frac{3}{2} \sqrt{(2 \pi)^{h+2n+2} n! (h+n+2)!} \frac{(-1)^j \beta_{h+n+2}}{(2 \pi)^j j! (n-j)! (h+n-j)!}.$$
From the definition of the inner product of the space $ \mathcal{A}_{n+1}(\mathbb{H})$ (see Definition \ref{inner}) we get
\begin{eqnarray*}
\langle B^{n+1}_{\mathcal{C}} \varphi, B^{n+1}_{\mathcal{C}} \psi \rangle_{\mathcal{A}_{n+1}(\mathbb{H})}&=& \sum_{h=0}^\infty \sum_{j=0}^n\frac{[(h+n-j)]^2(j!)^2((n-j)!)^2}{(2\pi)^{h+n+2-2j}(h+n+2)!} \overline{\beta_{h,j}'} \alpha_{h,j}' \\
&=&8 \sum_{h=0}^\infty \sum_{j=0}^n\frac{[(h+n-j)]^2(j!)^2((n-j)!)^2}{(2\pi)^{h+n+2-2j}(h+n+2)!}  \cdot\\
&& \cdot \frac{(2 \pi)^{h+2n+2} n! (h+n+2)!}{(2 \pi)^{2j}[j!(n-j)!]^2 [(h+n-j)]^2} \overline{\beta_{h+n+2}} \alpha_{h+n+2} \\
&=& 8 n! (2 \pi)^n  \sum_{j=0}^n\sum_{h=0}^\infty \overline{\beta_{h+n+2}} \alpha_{h+n+2} \\
&=& 8 n! (2 \pi)^n (n+1)\sum_{k=n+2}^\infty \overline{\beta}_k \alpha_k\\
&=& 8 (n+1)! (2 \pi)^n \langle \varphi, \psi \rangle_{ \mathcal{H}^n}.
\end{eqnarray*}
In particular when $ \varphi= \psi$ we obtain
$$ \| B^{n+1}_{\mathcal{C}} \varphi \|_{\mathcal{A}_{n+1}(\mathbb{H})}= 2 \sqrt{2} \sqrt{(n+1)!}(2 \pi)^{\frac{n}{2}} \| \varphi \|_{\mathcal{H}^n}.$$
\end{proof}

\begin{rem}
If we put $n=0$ in all the properties that we proved for $B^{n+1}_{\mathcal{C}}$: Proposition \ref{prop5}, Proposition \ref{ris3}, Proposition \ref{ris2}, respectively; we get up to a constant \cite[Prop. 4.18]{DKS2019}, \cite[Prop. 4.19]{DKS2019} and \cite[Prop. 4.20]{DKS2019}, respectively.
\end{rem}
We believe that it is also possible to define a $ \mathcal{C}$-polyanalytic Fueter Bargmann transform by applying the polyanalytic Fueter mapping $ \mathcal{C}_{n+1}$ to the reproducing kernel of the QTP Fock space, see \cite[Prop. 4.2]{DMD2}. In order to show how hard are the computations we show the case $n=1$. This case suggests that the formula that we are looking for puts in relation the quaternionic generalized Appell polyanalytic polynomials and the quaternionic Hermite polynomials.
\begin{prop}
\label{poly2}
For any $q,r \in \mathbb{H}$, we denote by $K_2(q,r)$ the reproducing kernel of the QTP Fock space, in the case $n=1$, then we have
$$
\mathcal{C}_2(K_2(q,r))= -8 \left( \sum_{h=0}^{\infty} \frac{\mathcal{M}_{1,h}(q, \bar{q}) \overline{H^{2 \pi}_{1, h+2}}(r, \bar{r})}{h!}-\sum_{h=0}^\infty \frac{\mathcal{M}_{0,h}(q, \bar{q})}{(2 \pi) h!} \overline{H_{1,h+3}^{2 \pi}(r, \bar{r})} \right).
$$
\end{prop}
\begin{proof}
Before to start we recall the definition of *-exponential
$$e_{*}(2 \pi q \bar{r})=\sum_{n=0}^\infty \frac{(2 \pi)^n q^n \bar{r}^n}{n!}$$
Then, by \cite[Prop. 4.2]{DMD2} with $n=1$ we have
\begin{eqnarray*}
K_2(q,r) &=& 2 e_{*}(2 \pi q \bar{r})* \left( \sum_{k=0}^1 (-1)^k \binom{1}{1-k} \frac{1}{k!}[2 \pi(q \bar{q}-q \bar{r}- \bar{q} r+ \bar{r} r)]^{k*}  \right)\\
&=& 2 e_{*}(2 \pi q \bar{r})* \left(1-2 \pi( \bar{q}q - q \bar{r}- \bar{q}r+ \bar{r} r) \right) \\
&=& 2 e_{*}(2 \pi q \bar{r})-4 \pi e_{*}(2 \pi q \bar{r}) |q|^2+4 \pi [ e_{*}(2 \pi q \bar{r})* q \bar{r}]+4 \pi [e_{*}(2 \pi q \bar{r})* \bar{q}r]+\\
&&-4 \pi e_{*}(2 \pi q \bar{r}) |r|^2.
\end{eqnarray*}
Now, we observe that
\begin{eqnarray*}
e_{*}(2 \pi q \bar{r})* q \bar{r} &=& 2 \left( \sum_{n=0}^\infty \frac{(2 \pi)^n q^n \bar{r}^n}{n!} \right)* q \bar{r}\\
&=& 2 \left( \sum_{n=0}^\infty \frac{(2 \pi)^n q^{n+1} \bar{r}^n}{n!} \right) \bar{r} \\
&=& q e_{*}(2 \pi q \bar{r}) \bar{r}.
\end{eqnarray*}
The same holds for the other member
$$ e_{*}(2 \pi q \bar{r})* \bar{q} r= \bar{q} e_{*}(2 \pi q \bar{r}) r.$$
Hence, we can write the reproducing kernel $K_2$ as the following polyanalytic decomposition
\begin{eqnarray}
\label{poly1}
K_2(q,r) &=&  2e_{*}(2 \pi q \bar{r})-4 \pi e_{*}(2 \pi q \bar{r})|r|^2+4 \pi q e_{*}(2 \pi q \bar{r}) \bar{r}\\
\nonumber
&& + \bar{q}[4 \pi e_{*}(2 \pi q \bar{r})r-4 \pi q e_{*}(2 \pi q \bar{r})].
\end{eqnarray}
Now, we apply the polyanalytic Fueter map $ \mathcal{C}_2$ to obtain
\begin{eqnarray}
\label{rel2}
\, \, \, \, \, \, \, \, \, \, \, \, \, \, \, \, \mathcal{C}_2(K_2(q,r)) &=& 2 \Delta \left(e_*(2 \pi q \bar{r})\right)- 4 \pi \Delta \left(e_*(2 \pi q \bar{r})| r|^2 \right)+4 \pi \Delta \left(q e_{*}(2 \pi q \bar{r})\right)\bar{r}\\ \nonumber
&& + 4 \pi x_0[\Delta \left(e_*(2 \pi q \bar{r}) r \right)-\Delta \left(q e_*(2 \pi q \bar{r})\right)].
\end{eqnarray}
By \cite[Prop. 4.2]{DKS2019}, up to a constant, we know that
\begin{equation}
\label{rel3}
\Delta\left(e_{*}(2 \pi q \bar{r}) \right)=-4 \sum_{h=0}^\infty \frac{(2 \pi)^{h+2} \mathcal{Q}_h(q, \bar{q}) \bar{r}^{h+2}}{h!}.
\end{equation}
Moreover,
\begin{eqnarray}
\label{rel4}
\Delta[q e_{*}(2 \pi q \bar{r})] &=& 2 \sum_{\ell=0}^\infty  \frac{(2 \pi)^{\ell} \Delta(q^{\ell+1}) \bar{r}^\ell}{\ell!} \\ \nonumber
&=& -4 \sum_{\ell=1}^{\infty}  \frac{(2 \pi)^{\ell} (\ell+1) \ell \mathcal{Q}_{\ell-1}(q, \bar{q}) \bar{r}^{\ell}}{ \ell (\ell-1)!} \\ \nonumber
&=& -4 \sum_{h=0}^\infty \frac{(2 \pi)^{h+1}(h+2)\mathcal{Q}_{h}(q, \bar{q}) \bar{r}^{h+1}}{h!}.
\end{eqnarray}
Now, we split formula \eqref{rel2} in two parts
$$ A:=2\Delta \left(e_*(2 \pi q \bar{r})\right)- 4 \pi \Delta \left(e_*(2 \pi q \bar{r})| r|^2 \right)+4 \pi \Delta \left(q e_{*}(2 \pi q \bar{r})\right) \bar{r},$$
$$ B:= 4 \pi x_0[\Delta \left(e_*(2 \pi q \bar{r}) r \right)-\Delta \left(q e_*(2 \pi q \bar{r})\right)],$$
thus
$$ \mathcal{C}_2(K_2)(q,r)=A+B.$$
Firstly, we calculate $B$. By formula \eqref{rel3} and formula \eqref{rel4} we get
\begin{eqnarray*}
B &=& 4 \pi x_0 \left(-4 \sum_{h=0}^\infty \frac{(2 \pi)^{h+2} \mathcal{Q}_h(q, \bar{q}) \bar{r}^{h+2}}{h!} \cdot r+4 \sum_{h=0}^\infty  \frac{(2 \pi)^{h+1}(h+2) \mathcal{Q}_h(q, \bar{q}) \bar{r}^{h+1}}{h!} \right) \\
&=& -8 \left[ \sum_{h=0}^\infty \frac{x_0 \mathcal{Q}_h(q, \bar{q})}{h!} \left((2 \pi)^{h+3} r \bar{r}^{h+2}-(2 \pi)^{h+2}(h+2) \bar{r}^{h+1} \right) \right] \\
&=& -8 \sum_{h=0}^\infty \frac{\mathcal{M}_{1,h}(q, \bar{q})}{h!} \overline{H^{2 \pi}_{1,h+2}(r, \bar{r})}.
\end{eqnarray*}
Now, we calculate $A$.
\begin{eqnarray*}
A &=& 2 \biggl(-4 \sum_{h=0}^\infty \frac{(2 \pi)^{h+2} \mathcal{Q}_h(q, \bar{q}) \bar{r}^{h+2}}{h!}+8 \pi \sum_{h=0}^\infty \frac{(2 \pi)^{h+2} \mathcal{Q}_h(q, \bar{q}) \bar{r}^{h+2} |r |^2}{h!}+\\
&& -8 \pi \sum_{h=0}^{\infty} \frac{(2 \pi)^{h+1} (h+2) \mathcal{Q}_h(q, \bar{q}) \bar{r}^{h+2}}{h!} \biggl)\\
&=& 8 \left[ \sum_{h=0}^\infty \frac{\mathcal{Q}_h(q, \bar{q})}{(2 \pi) h!} \left( (2 \pi)^{h+4} r \bar{r}^{h+3}-(2 \pi)^{h+3}(h+2+1) \bar{r}^{h+2} \right) \right]\\
&=& 8 \sum_{h=0}^\infty \frac{\mathcal{M}_{0,h}(q, \bar{q})}{(2 \pi) h!} \overline{H_{1,h+3}^{2 \pi}(r, \bar{r})}.
\end{eqnarray*}
Hence
\begin{eqnarray*}
\mathcal{C}_2(K_2(q,r)) &=& A+B \\
&=& -8 \left( \sum_{h=0}^{\infty} \frac{\mathcal{M}_{1,h}(q, \bar{q}) \overline{H^{2 \pi}_{1, h+2}(r, \bar{r})}}{h!}-\sum_{h=0}^\infty \frac{\mathcal{M}_{0,h}(q, \bar{q})}{(2 \pi) h!} \overline{H_{1,h+3}^{2 \pi}(r, \bar{r})} \right).
\end{eqnarray*}
\end{proof}
\begin{rem}
Since the computations are already hard for the case $n=1$, we do not know if it is possible to find a formula for all $n$. This would be investigated in a forthcoming work.
\end{rem}
\subsection{The map $ \tau_{n+1}$ applied to the QTP Bargmann transform}
In this subsection we want to investigate what happens if we apply the polyanalytic Fueter mapping $\tau_{n+1}$ to the QTP Bargmann $B^{n+1}$. We will call this integral transforms $\tau$- polyanalytic Fueter-Bargmann transform.
\\ First of all we want to study the range of this integral transform. Let us define the following space
\begin{equation}
\label{for2}
\widetilde{\mathcal{A}}_{n+1}(\mathbb{H}):= \left \{ \tau_{n+1}(f); \, f \in \mathcal{F}_n^T(\mathbb{H}) \right \}.
\end{equation}
Before to give a characterization of this space. We need the following result
\begin{prop}
\label{apph}
Let $n\geq 0$ a fixed number. Then we have
$$ \tau_{n+1}\left(H^{2 \pi}_{n,j}(q, \bar{q})\right)=
\begin{cases}
-2^{n+1}(2 \pi)^{j+n} n! j(j-1) \mathcal{Q}_{j-2}(q, \bar{q})& \, \, \, \, j\geq n+2, \\
0 &  \, \, \, \, j <n+2.
\end{cases}
$$ .
\end{prop}
\begin{proof}
By using the relation between the two polyanalytic Fueter maps, see \cite[Thm. 3.13]{ADS} and Proposition \ref{prop2} we have for $j \geq n+2$
\begin{eqnarray}
\label{for1}
\nonumber
\tau_{n+1}\left(H_{n,j}^{2 \pi}(q, \bar{q}) \right) &=& 2^{n} \mathcal{D}^{n} \mathcal{C}_{n+1}\left( H_{n,j}^{2\pi}(q, \bar{q}) \right) \\
&=& -2 ^{n}2 (2 \pi)^{j} n! j! \sum_{s=0}^n \frac{(-1)^s (2 \pi)^{n-s}}{s! (n-s)!(j-s-2)!}\mathcal{D}^{n} \left( \mathcal{M}_{n-s, j-s-2}(q, \bar{q}) \right).
\end{eqnarray}
From Proposition \ref{res1} we know that
$$ \mathcal{D}^n  \left( \mathcal{M}_{n-s, j-s-2}(q, \bar{q}) \right)= n! \mathcal{Q}_{j-s-2}(q, \bar{q}),$$
if $n-s=n$, this means $s=0$. Therefore
$$ \mathcal{D}^{n}\left( \mathcal{M}_{n-s, j-s-2}(q, \bar{q}) \right)= n! \mathcal{Q}_{j-2}(q, \bar{q}).$$
This implies that in the sum \eqref{for1} we have only to take into account the case $s=0$. Therefore for $j \geq n+2$ we get
\begin{eqnarray*}
\tau_{n+1}\left(H_{n,j}^{2 \pi}(q, \bar{q}) \right)&=& - \frac{2^n 2 (2 \pi)^{j+n}(n!)^2 j!}{n! (j-2)!} \mathcal{Q}_{j-2}(q, \bar{q}) \\
&=& -2 ^{n+1}(2 \pi)^{j+n} n!j(j-1) \mathcal{Q}_{j-2}(q, \bar{q}).
\end{eqnarray*}
\end{proof}
Now, we can characterize the space $ \widetilde{\mathcal{A}}_{n+1}(\mathbb{H})$.
\begin{thm}
\label{thm24}
For any fixed $n\geq 0$, we have
$$ \widetilde{\mathcal{A}}_{n+1}(\mathbb{H})= \left \{ \sum_{h=0}^\infty \mathcal{Q}_{h+n}(q, \bar{q}) \beta_h, \, \, ( \beta_h)_{h \geq 0} \subset \mathbb{H}, \, \,  \sum_{h=0}^\infty \frac{(h+n)!| \beta_h|^2}{(h+n+1)(h+n+2)(2 \pi)^h}< \infty \right \}.$$
\end{thm}
\begin{proof}
Let us assume that $g \in \widetilde{\mathcal{A}}_{n+1}(\mathbb{H})$, thus there exists a function $f \in \mathcal{F}_n^T(\mathbb{H})$ such that $g(q)= \tau_{n+1}(f)(q)$. Since $f \in \mathcal{F}_n^T(\mathbb{H})$ by Proposition \ref{prop1} we have
$$ f(q)= \sum_{j=0}^\infty H_{n,j}(q, \bar{q})^{2 \pi} \alpha_j, \qquad \{\alpha_{j}\}_{j\geq 0} \subset \mathbb{H}$$
and
$$ \| f \|_{\mathcal{F}_n^T(\mathbb{H})}^2= \sum_{j=0}^{\infty} (2 \pi)^j j! | \alpha_j|^2 < \infty.$$
Now, by Proposition \ref{apph} we get
\begin{eqnarray*}
g(q)&=& \tau_{n+1}(f)(q)= \sum_{j=0}^\infty \tau_{n+1} \left(H_{n,j}(q, \bar{q}) \right) \alpha_j \\
&=& -\sum_{j=n+2}^\infty 2^{n+1} (2 \pi)^{j+n} n! j(j-1) \mathcal{Q}_{j-2}(q, \bar{q}) \alpha_j \\
&=& -2^{n+1} (2 \pi)^{2(n+1)} n! \sum_{h=0}^{\infty} (2 \pi)^{h} (h+n+2)(h+n+1) \mathcal{Q}_{h+n}(q, \bar{q}) \alpha_{h+n+2} \\
&=& \sum_{h=0}^\infty \mathcal{Q}_{h+n}(q, \bar{q}) \beta_h,
\end{eqnarray*}
where $ \beta_h:=-2^{n+1} (2 \pi)^{2(n+1)} n! (h+n+2)(h+n+1) (2 \pi)^h \alpha_{h+n+2}$.
Moreover
\begin{eqnarray*}
\sum_{h=0}^\infty \frac{(h+n)!}{(h+n+1)(h+n+2)(2 \pi)^h}| \beta_h|^2 &=& \sum_{h=0}^\infty \frac{(h+n)!}{(h+n+1)(h+n+2)(2 \pi)^h} 4^{n+1}(2 \pi)^{4(n+1)}(n!)^2 \cdot \\
&& \cdot (h+n+1)^2(h+n+2)^2 (2 \pi)^{2h} | \alpha_{h+n+2}|^2 \\
&=&  4^{n+1}(2 \pi)^{4(n+1)}(n!)^2 \sum_{h=0}^\infty (h+n)! (h+n+1)(h+n+2) \cdot\\
&& \cdot (2 \pi)^h | \alpha_{h+n+2}|^2\\
&=& 4^{n+1}(2 \pi)^{3n+2}(n!)^2  \sum_{h=0}^\infty (h+n+2)! (2 \pi)^{h+2+n} | \alpha_{h+n+2}|^2\\
&=& 4^{n+1}(2 \pi)^{3n+2}(n!)^2 \| f \|_{\mathcal{F}_n^T(\mathbb{H})}^2 < \infty.
\end{eqnarray*}
Therefore
$$ \sum_{h=0}^\infty \frac{(h+n)!}{(h+n+1)(h+n+2)(2 \pi)^h}| \beta_h|^2 < \infty. $$
Now, we prove the other inclusion. Let us consider a function
$$ h(q)= \sum_{k=n+2}^\infty H_{n,k}^{2 \pi} (q, \bar{q}) \alpha_k,$$
with $ \alpha_k:=- \frac{\beta_{k-n-2}}{2^{n+1} (2 \pi)^{n+k} n! (k-1)k}.$
We get that $g(q)= \tau_{n+1}(h)(q)$, since
$$ \mathcal{Q}_{k+n}(q, \bar{q})=- \frac{\tau_{n+1} \left(H_{n,k+n+2}^{2\pi}(q, \bar{q}) \right)}{2^{n+1}(2 \pi)^{k+2n+2}n! (k+n+2)(k+n+1)}.$$
Furthermore,
\begin{eqnarray*}
\| h \|_{\mathcal{F}_n^T(\mathbb{H})}^2 &=& \sum_{k=n+2}^\infty (2 \pi)^k k! | \alpha_k|^2 \\
&=& \sum_{\ell=0}^\infty (2 \pi)^{\ell+n+2} ( \ell+n+2)! | \alpha_{\ell+n+2}|^2 \\
&=& \sum_{\ell=0}^\infty  \frac{(2 \pi)^{\ell+n+2} ( \ell+n+2)!| \beta_{\ell}|^2}{4^{n+1}(2 \pi)^{2(\ell+2n+2)} (n!)^2 ( \ell+n+1)^2 ( \ell+n+2)^2} \\
&=& \frac{1}{4^{n+1} (n!)^2 (2 \pi)^{3n+2}} \sum_{\ell=0}^\infty \frac{(\ell+n)!}{( \ell+n+1)(\ell+n+2) (2 \pi)^ \ell} | \beta_{\ell}| ^2 < \infty.
\end{eqnarray*}
Hence $h \in \mathcal{F}_n^T(\mathbb{H})$ and $g(q)= \tau_{n+1}(h)(q)$. This means that $g \in \widetilde{A}_{n+1} (\mathbb{H})$.
\end{proof}
\begin{defn}
\label{ad1}
Let $n \geq 0$ be a fixed number. Let us consider $f(q)=\sum_{h=0}^\infty \mathcal{Q}_{h+n}(q, \bar{q}) \beta_h$ and $g(q)=\sum_{h=0}^\infty \mathcal{Q}_{h+n}(q, \bar{q}) \gamma_h$ we define the inner product of $\widetilde{\mathcal{A}}_{n+1} (\mathbb{H})$ as
$$ \langle f,g \rangle_{\widetilde{A}_{n+1} (\mathbb{H})}= \sum_{k=0}^\infty \frac{(k+n)!}{( k+n+1)(k+n+2) (2 \pi)^ k} \overline{\gamma_k} \beta_k.$$
The associated norm is given by
$$ \|f\|^2_{\widetilde{\mathcal{A}}_{n+1} (\mathbb{H})}=\sum_{k=0}^\infty \frac{(k+n)!}{( k+n+1)(k+n+2) (2 \pi)^ k} |\beta_k|^2.$$
\end{defn}
Now, we can propose the following
\begin{defn}[$\tau$- polyanalytic Fueter Bargmann transform]
\label{tb}
We define $ B^{n+1}_{\tau}: L^{2}_{\mathbb{H}}( \mathbb{R}) \mapsto \widetilde{A}_{n+1}( \mathbb{H})$ as
$$ B^{n+1}_{\tau}:= \tau_{n+1} \circ B^{n+1},$$
where $B^{n+1}$ is the QTP Bargmann transform defined in formula \eqref{pre}.
\end{defn}

\begin{cor}
The $\tau$- polyanalytic Fueter-Bargmann transform $B^{n+1}_{\tau}$ is $ \mathbb{H}$-linear.
\end{cor}
\begin{proof}
This result follows by the linearity of the map $\tau_{n+1}$, see Proposition \ref{line2}.
\end{proof}
It is possible to prove a relation between the $\tau$- polyanalytic Fueter-Bargmann transform and the $ \mathcal{C}$- polyanalytic Fueter-Bargmann transform.
\begin{prop}
\label{con}
Let $ \varphi \in L^{2}_{\mathbb{H}}( \mathbb{R})$. Then, for any fixed $n\geq 0$ we have
$$ B^{n+1}_{\tau} \varphi=2^{n} \mathcal{D}^n \left(B^{n+1}_{\mathcal{C}} \varphi \right).$$
\end{prop}
\begin{proof}
The result follows from \cite[Thm. 3.13]{ADS}, Definition \ref{tb} and Definition \ref{CB}
$$ B^{n+1}_{\tau} \varphi= \tau_{n+1} \circ (B^{n+1} \varphi)= 2^n \mathcal{D}^n \left( \mathcal{C}_{n+1} \circ B^{n+1}(\varphi) \right)=2^{n} \mathcal{D}^{n}(B^{n+1}_{\mathcal{C}} \varphi).$$
\end{proof}
As the $ \mathcal{C}$- polyanalytic Fueter-Bargmann transform we can express the $\mathcal{\tau}$- polyanalytic Fueter-Bargmann tranform in integral form.
\begin{thm}
For any fixed $n\geq 0$ it is possible to define the $\mathcal{\tau}$- polyanalytic Fueter-Bargmann through this diagram
$$B_{\tau}^{n+1}: \xymatrix{
	L^2_\Hq(\R) \ar[r] \ar[d]_{B^{n+1}} & \mathcal{\tilde{A}}_{n+1}(\Hq)  \\ \mathcal{F}^{n}_{T}(\Hq) \ar[r]_{Id} & \mathcal{SP}_{n+1}(\Hq) \ar[u]_{\tau_{n+1}}
}$$

Precisely, for any $ \varphi \in L^{2}_{\mathbb{H}}(\mathbb{R})$ and $q \in \mathbb{H}$ we have
$$ B^{n+1}_{\tau} \varphi(q)= \int_{\mathbb{R}} \Theta(q,x) \varphi(x) \, dx,$$
where
\begin{equation}
\label{tau1}
\Theta(q,x)=- \sqrt{n!} \pi^{\frac{n}{2}} 2^{\frac{3}{4}} \sum_{k=0}^\infty \frac{ \mathcal{Q}_{k+n}(q, \bar{q}) h^{2 \pi}_{k+n+2}(x)}{2^{\frac{k}{2}}(k+n)!}.
\end{equation}
\end{thm}
\begin{proof}
By Proposition \ref{con} we have that
$$ B^{n+1}_{\tau} \varphi(q)=2^{n} \mathcal{D}^n \left(B^{n+1}_{\mathcal{C}} \varphi(q) \right).$$
From Theorem \ref{CB1} we have that
$$ B^{n+1}_{\tau} \varphi (q)= 2^{n} \int_{\mathbb{R}} \mathcal{D}^n \left( \Phi(q,x) \right) \varphi(x) dx, \quad q\in\mathbb{H}.$$
This means that in order to prove statement it is enough to compute $\mathcal{D}^n \left( \Phi(q,x) \right)$, where
$$\Phi(q,x):=-\sqrt{n!}(2\pi)^{\frac{n}{2}}2^{\frac{3}{4}}\sum_{l=0}^{\infty}\sum_{j=0}^{n}\frac{(-1)^jh_{l+n+2}^{2\pi}(x)\mathcal{M}_{n-j,l+n-j}(q,\bar{q})}{2^{\frac{l+n}{2}}(2\pi)^j j!(l+n-j)!(n-j)!}.$$
Firstly, we observe that by Proposition \ref{res1}, for $k >0$, we have that
$$ \mathcal{D}^n ( \mathcal{M}_{n-j,k}(q, \bar{q}))=
\begin{cases}
n! \mathcal{Q}_{k}(q, \bar{q}) \qquad \hbox{if} \, \,j=0,\\
0 \qquad \qquad \qquad \hbox{if} \, \, j>0.
\end{cases}$$
This imply that when we apply the the operator $ \mathcal{D}^n$ to $ \Phi(q,x)$ the only term that we have to take into account in the second series is when $j=0$, since the others members are all zero. Now, we perform the computation
\begin{eqnarray*}
\mathcal{D}^n \left(\Phi(q,x)\right)&=&-\sqrt{n!}(2\pi)^{\frac{n}{2}}2^{\frac{3}{4}}\sum_{k=0}^{\infty}\sum_{j=0}^{n}\frac{(-1)^j h_{k+n+2}^{2\pi}(x)\mathcal{D}^n\left(\mathcal{M}_{n-j,k+n-j}(q,\bar{q})\right)}{2^{\frac{k+n}{2}}(2\pi)^j j!(k+n-j)!(n-j)!}\\
&=& - \sqrt{n!} \pi^{\frac{n}{2}} 2^{\frac{3}{4}} \sum_{k=0}^\infty \frac{ \mathcal{Q}_{k+n}(q, \bar{q}) h^{2 \pi}_{k+n+2}(x)}{2^{\frac{k}{2}}(k+n)!}.
\end{eqnarray*}
\end{proof}
\begin{rem}
If we put $n=0$ in \eqref{tau1} we obtain the same kernel of the integral transform obtained in \cite[Thm. 4.14]{DKS2019}, up to a constant.
\end{rem}
Now, we prove some properties of the $\tau$- poly Fueter-Bargmann transform .
\\ First of all we recall the normalized Hermite functions
$$ \psi_k^{2 \pi}(x)= \frac{h_{k}^{2 \pi}(x)}{ \| h_k^{2 \pi} \|_{L^{2}_{\mathbb{H}}(\mathbb{R})}},$$
where $h_k^{2 \pi}$ are defined in formula \eqref{hermi}.
\begin{prop}
\label{tp1}
For any fixed $n\geq 0$, the action of the $\tau$- poly Fueter-Bargmann transform on the normalized Hermite function is the following
$$ B^{n+1}_{\tau}(\psi_k^{2 \pi})(q)=
\begin{cases}
-\sqrt{2}2^{n+1}\sqrt{(2 \pi)^{k+n} n! k!} \frac{\mathcal{Q}_{k-2}(q, \bar{q})}{(k-2)!}  & \mbox{if } \quad k \geq n+2 \\
0 & \mbox{if} \quad k<n+2.
\end{cases}
$$
\end{prop}
\begin{proof}
From Proposition \ref{con} we know that
$$ B^{n+1}_{\tau}(\psi_k^{2 \pi})(q)=2^n \mathcal{D}^n \left(B^{n+1}_{\mathcal{C}} \psi_k^{2 \pi} \right)(q).$$
Now, if $k \geq n+2$, by Proposition \ref{prop5} we get that
$$ B^{n+1}_{\tau}( \psi_{k}^{2 \pi})(q)=-\sqrt{2}2^{n+1} \sqrt{(2 \pi)^{k+n} k! n!} \sum_{j=0}^n \frac{(-1)^j \mathcal{D}^n \left(\mathcal{M}_{n-j,k-j-2}(q, \bar{q}) \right)}{(2 \pi)^j j! (n-j)! (k-j-2)!}.$$
By Proposition \ref{res1} we have that
$$ \mathcal{D}^n \left(\mathcal{M}_{n-j,k-j-2}(q, \bar{q}) \right)=
\begin{cases}
n! \mathcal{Q}_{k-2}(q, \bar{q})& \qquad \hbox{if} \, \, j=0,\\
0& \qquad   \, \, \hbox{if} \, \, j> 0.
\end{cases}$$
This means that in the summation survives only the term with $j=0$. Hence, we get
\begin{eqnarray*}
B^{n+1}_{\tau}(\psi_k^{2 \pi})(q) &=& -\sqrt{2}2^{n+1} \sqrt{(2 \pi)^{k+n} k! n!} \frac{n! \mathcal{Q}_{k-2}(q, \bar{q})}{n! (k-2)!} \\
&=& -\sqrt{2}2^{n+1} \sqrt{(2 \pi)^{k+n} k! n!} \frac{ \mathcal{Q}_{k-2}(q, \bar{q})}{ (k-2)!}.
\end{eqnarray*}
\end{proof}
\begin{rem}
If we put $n=0$ in Proposition \ref{tp1} we get the same result of \cite[Prop. 4.18]{DKS2019}, up to a constant. Indeed, if $k \geq 2$ we get
\begin{eqnarray*}
B^{1}_{\tau}( \psi_k^{2 \pi})(q) &=& -2 \sqrt{2} \sqrt{(2 \pi)^k k!} \frac{\mathcal{Q}_{k-2}(q, \bar{q})}{(k-2)!}\\
&=& -2 \sqrt{2} \sqrt{(2 \pi)^k k(k-1)} \sqrt{(k-2)!} \frac{\mathcal{Q}_{k-2}(q, \bar{q})}{(k-2)!}\\
&=& -2 \sqrt{2} \sqrt{\frac{(2 \pi)^k k (k-1)}{(k-2)!}} \mathcal{Q}_{k-2}(q, \bar{q}).
\end{eqnarray*}
\end{rem}
Now, using Definition \ref{hil}, we have the following

\begin{prop}
\label{tp2}
Let $ \varphi\in \mathcal{H}^n$. Then, for any fixed $n\geq 0$ we have
$$B^{n+1}_{\tau} (\varphi)(q)=\sum_{k=0}^\infty \mathcal{Q}_{k+n}(q, \bar{q}) \beta_k,$$
where
$$ \beta_k=- \frac{\sqrt{2}2^{n+1} \sqrt{(2 \pi)^{k+2+2n} n! (k+n+2)!}}{(k+n)!} \alpha_{k+n+2}, \qquad \{ \alpha_k \}_{k\geq 0 } \subset \mathbb{H}.$$
\end{prop}
\begin{proof}
By Definition \ref{hil} we can write
$$ \varphi(x)= \sum_{k = n+2}^\infty \psi_k^{2 \pi}(x) \alpha_k, \qquad \{ \alpha_k \}_{k\geq 0} \subset \mathbb{H}.$$
Then by Proposition \ref{tp1} we get
\begin{eqnarray*}
B^{n+1}_{\tau}( \varphi)(q) &=& \sum_{k=n+2}^\infty B^{n+1}_{\tau} (\psi_k^{2 \pi})(q) \alpha_k \\
&=& - \sum_{k=n+2}^\infty  \sqrt{2}2^{n+1} \sqrt{(2 \pi)^{k+n} n! k!} \frac{\mathcal{Q}_{k-2}(q, \bar{q})}{(k-2)!} \alpha_k \\
&=& - \sum_{k=0}^\infty \sqrt{2}2^{n+1} \sqrt{(2 \pi)^{k+2+2n} n! (k+n+2)!} \frac{\mathcal{Q}_{k+n}(q, \bar{q})}{(k+n)!} \alpha_{k+n+2}  \\
&:=& \sum_{k=0}^\infty \mathcal{Q}_k(q, \bar{q}) \beta_k.
\end{eqnarray*}
\end{proof}

\begin{thm}
\label{tp4}
The $\tau$-polyanalytic Bargmann transform is a surjective operator that satisfies the following unitary and isometric properties for $ \varphi , \psi \in \mathcal{H}^n$
$$ \langle B^{n+1}_{\tau} \varphi, B^{n+1}_{\tau} \psi\rangle_{\widetilde{\mathcal{A}}_{n+1}(\mathbb{H})}=4^{n+2} 2(2 \pi)^{2(n+1)} n! \langle \varphi, \psi \rangle_{\mathcal{H}^n}.$$
In particular
\begin{equation}
\label{tp3}
\| B^{n+1}_{\tau} \varphi \|_{\mathcal{\tilde{A}}_{n+1}(\mathbb{H})}= \sqrt{2}2^{n+2} (2 \pi)^{n+1} \sqrt{n!} \| \varphi \|_{\mathcal{H}^n}.
\end{equation}
\end{thm}
\begin{proof}
Let us consider
$$ \varphi=\sum_{k=n+2}^\infty \psi_k^{2 \pi}\alpha_k, \qquad \psi=\sum_{k=n+2}^\infty \psi_k^{2 \pi} \alpha_k',$$
with $ \{\alpha_k\}_{k\geq 0} \subset \mathbb{H}$ and $ \{ \alpha_k'\}_{k\geq 0}\subset \mathbb{H}$. From Proposition \ref{tp2} we know that
$$B^{n+1}_{\tau} (\varphi)(q)= \sum_{k=0}^\infty \mathcal{Q}_{k+n}(q, \bar{q}) \beta_k\qquad \hbox{and} \qquad B^{n+1}_{\tau} (\psi)(q)= \sum_{k=0}^\infty \mathcal{Q}_{k+n}(q, \bar{q}) \gamma_k,$$
where
$$ \beta_k:= - \frac{\sqrt{2}2^{n+1} \sqrt{(2 \pi)^{k+2+2n} n! (k+n+2)!}}{(k+n)!} \alpha_{k+n+2}, $$
$$\gamma_k:=- \frac{\sqrt{2}2^{n+1} \sqrt{(2 \pi)^{k+2+2n} n! (k+n+2)!}}{(k+n)!} \alpha'_{k+n+2}.$$
Now, recalling the definition of the inner product $ \widetilde{\mathcal{A}}_{n+1}(\mathbb{H})$ (see Definition \ref{ad1}) we get
\begin{eqnarray*}
\langle B^{n+1}_{\tau} \varphi, B^{n+1}_{\tau} \psi\rangle_{\widetilde{\mathcal{A}}_{n+1}(\mathbb{H})}&=& \sum_{k=0}^\infty \frac{(k+n)!}{(k+n+2)(k+n+1)(2 \pi)^k} \overline{\gamma}_k \beta_k\\
&=& \sum_{k=0}^\infty \frac{(k+n)!}{(k+n+2)(k+n+1)(2 \pi)^k}  \frac{4^{n+2}2 (2 \pi)^{2+k+2n} n! (k+n+2)!}{[(k+n)!]^2} \cdot\\
&&\cdot \overline{\alpha'_{k+n+2}} \alpha_{k+n+2}\\
&=&4^{n+2} 2(2 \pi)^{2(n+1)} n!\sum_{k=0}^\infty \frac{(k+n)!}{(k+n+2)(k+n+1)(2 \pi)^k} \cdot    \\
&& \cdot \frac{(2 \pi)^k (k+n+2)(k+n+1)(k+n)!}{[(k+n)!]^2}\overline{\alpha'_{k+n+2}} \alpha_{k+n+2}\\
&=& 4^{n+2} 2(2 \pi)^{2(n+1)} n! \sum_{k=n+2}^\infty \overline{\alpha_k'}\alpha_k\\
&=&4^{n+2} 2(2 \pi)^{2(n+1)} n!  \langle \varphi, \psi \rangle_{\mathcal{H}^n}.
\end{eqnarray*}
In particular if $\psi= \varphi$ we get formula \eqref{tp3}.
\end{proof}
\begin{cor}
\label{tp5}
The $\tau$-poly Fueter-Bargmann transform is a quaternionic bounded operator such that for any $ \varphi \in L^2_{\mathbb{H}}(\mathbb{R})$, we have
$$ \| B^{n+1}_{\tau} \|_{\mathcal{\tilde{A}}_{n+1}(\mathbb{H})} \leq \sqrt{2}2^{n+2} (2 \pi)^{n+1} \sqrt{n!} \| \varphi \|_{L^2_{\mathbb{H}}(\mathbb{R})}.$$
\end{cor}
\begin{proof}
Since $ \mathcal{H}^n$ is a subspace of $L^{2}_{\mathbb{H}}(\mathbb{R})$ we get that
$$ \| \varphi \|_{\mathcal{H}^n} \leq \| \varphi \|_{L^{2}_{\mathbb{H}}(\mathbb{R})}.$$
Therefore by applying this inequality to \eqref{tp3} we get the thesis.
\end{proof}
\begin{rem}
If we put $n=0$ in Proposition \ref{tp1}, Theorem \ref{tp4} and Corollary \ref{tp5} we get \cite[Prop. 4.18]{DKS2019}, \cite[Prop. 4.20]{DKS2019} and \cite[Prop. 4.19]{DKS2019}, respectively.
\end{rem}
Also in this case it is possible to construct a $\tau$-polyanalytic Fueter Bargmann transform by applying the Fueter poly map $\tau_{n+1}$ to the quaternionic reproducing kernel of the QTP Fock space $K_{n+1}$. However, also in this case the computations are very hard. For this reason we write only the case $n=1$.
\begin{prop}
For any $q,r \in \mathbb{H}$ we denote by $K_2(q,r)$ the reproducing kernel of the QTP Fock space, in the case $n=1$, then we have
\begin{equation}
\label{new}
\tau_{2} \left(K_{2}(q, r)\right)=-8 \sum_{k=0}^\infty \frac{\mathcal{Q}_{k}(q)}{k!} \overline{H^{2 \pi}_{1, k+2}(r, \bar{r})}.
\end{equation}
\end{prop}
\begin{proof}
From formula \eqref{poly1} we have the polyanalytic decomposition of $K_{2}(q,r)$. Recalling the action of the global operator $V$, we get that
$$ V \left( K_{2}(q, r)\right)=4 \pi [e_{*}(2 \pi q \bar{r})r-q e_{*}(2 \pi q \bar{r})].$$
Using similar computations as Proposition \ref{poly2} we get
\begin{eqnarray*}
\tau_{2} \left(K_{2}(q, r)\right)&=& \Delta\circ V \left(K_{2}(q,r)\right)\\
&=&4 \pi [\Delta \left( e_{*}(2 \pi q \bar{r})r\right)-\Delta \left(q e_{*}(2 \pi q \bar{r})\right)]\\
&=& -8 \sum_{k=0}^\infty \frac{\mathcal{Q}_{k}(q)}{k!} \overline{H^{2 \pi}_{1, k+2}(r, \bar{r})}.
\end{eqnarray*}
\end{proof}

\begin{rem}
It will be interesting to find a general expression for the formula \eqref{new} in the case $ n \geq 2$.
\end{rem}

\section{The polyanalytic Fueter Bargmann transforms}
In this section we see how the results obtained for the QTP Fock space can be reformulated for the quaternionic polyanalytic Fock space. In order to do this it is crucial the following result, see \cite[Thm. 3.4]{DMD2},
\begin{thm}
\label{sum}
Let $N \geq 0$. The quaternionic polyanalytic Fock space $ \mathcal{\widetilde{F}}_{Slice}^{N+1}(\mathbb{H})$ is the direct sum of QTP Fock spaces $\mathcal{F}_T^n(\mathbb{H})$, $n=0,...,N$ i.e.
$$ \mathcal{\widetilde{F}}_{Slice}^{N+1}(\mathbb{H})= \bigoplus_{n=0}^N   \mathcal{F}_T^n(\mathbb{H}).$$
\end{thm}
\begin{rem}
From the Theorem \ref{sum} it is clear that a function $f \in \mathcal{\widetilde{F}}_{Slice}^{N+1}(\mathbb{H})$ if and only if
$$ f(q)=\sum_{n=0}^N f_{n}(q), \quad f_n \in \mathcal{F}_T^n(\mathbb{H}), \, n=0,...,N.$$
\end{rem}
We will omit all proofs because by the previous remark they are similar to those ones obtained in the previous section.
\\ First of all we characterize the space $\mathcal{\widetilde{F}}_{Slice}^{N+1}(\mathbb{H})$. By putting the sum from $n=0$ to $N$ in the proof of Proposition \ref{prop1} we get the following
\begin{prop}
A function of the form
$$ f(q)= \sum_{n=0}^N \sum_{j=0}^\infty H_{n,j}^{2 \pi}(q, \bar{q}) \alpha_{n,j} \quad \{\alpha_{n,j}\}_{0 \leq n \leq N, j\geq 0} \subset \mathbb{H}, $$ belongs to the space $\mathcal{\widetilde{F}}_{Slice}^{N+1}(\mathbb{H})$ if and only if
$$ \sum_{n=0}^N \sum_{j=0}^\infty (2 \pi)^{j+n} j! n! | \alpha_{n,j}|^2 < \infty.$$
\end{prop}
\subsection{The map $ \mathcal{C}_{n+1}$ applied to the QFP Bargmann transform}
Let us start by defining the range of the polyanalytic Fueter mapping $ \mathcal{C}_{n+1}$ on the quaternionic polyanalytic Fock space,
$$ \mathfrak{A}_{N}(\mathbb{H})= \bigoplus_{n=0}^N \mathcal{A}_{n+1}(\mathbb{H}),$$
where $\mathcal{A}_{n+1}(\mathbb{H})$ is the space defined in formula \eqref{ff}.
\\We have the following characterization of the previous space
\begin{thm}
Let $N \geq 0$. Then we have
\begin{eqnarray*}
\mathfrak{A}_{N}(\mathbb{H})&=& \! \!  \biggl \{ \sum_{n=0}^N \sum_{h=0}^\infty \sum_{s=0}^n \mathcal{M}_{n-s,h+n-s}(q,\bar{q}) \beta_{n,h,s},\\
&& \sum_{n=0}^N \sum_{h=0}^\infty \sum_{s=0}^n \frac{[(h+n-s)]^2(s!)^2[(n-s)!]^2}{(n+1)!(2\pi)^{h+2n+2-2s}(h+n+2)!}| \beta_{n,h,s}|^2 < \infty\biggl \}.
\end{eqnarray*}

where $\{\beta_{n,h,s}\}_{0 \leq n \leq N, h\geq 0, 0 \leq s \leq n } \subset \mathbb{H}$.
\end{thm}

\begin{defn}
Let $f,g \in \mathfrak{A}_{N}(\mathbb{H})$ be such that $f=\sum_{n=0}^N \sum_{h=0}^\infty \sum_{s=0}^n \mathcal{M}_{n-s,h+n-s}(q, \bar{q}) \alpha_{n,h,s}$ and $g=\sum_{n=0}^N \sum_{h=0}^\infty \sum_{s=0}^n \mathcal{M}_{n-s,h+n-s}(q, \bar{q}) \beta_{n,h,s}$. We define their inner product as
$$ \langle f,g \rangle_{\mathfrak{A}_{N}(\mathbb{H})}= \sum_{n=0}^N \frac{\langle f,g \rangle_{\mathcal{A}_{n+1}(\mathbb{H})}}{(n+1)! (2 \pi)^{n}}= \sum_{n=0}^N \sum_{h=0}^\infty \sum_{s=0}^n \frac{[(h+n-s)]^2(s!)^2[(n-s)!]^2}{(n+1)!(2\pi)^{h+2n+2-2s}(h+n+2)!}\overline{ \beta_{n,h,s}}\alpha_{n,h,s} $$
and the norm as
$$ \| f \|_{\mathfrak{A}_{N}(\mathbb{H})}^2= \sum_{n=0}^N \frac{\|f \|_{\mathcal{A}_{n+1}(\mathbb{H})}^2}{(n+1)! (2 \pi)^{n}}=\sum_{n=0}^N \sum_{h=0}^\infty \sum_{s=0}^n \frac{[(h+n-s)]^2(s!)^2[(n-s)!]^2}{(n+1)!(2\pi)^{h+2n+2-2s}(h+n+2)!} |\alpha_{n,h,s}|^2.$$
\end{defn}
Now, we can give the following
\begin{defn}[$ \mathcal{C}$- full poly Fueter Bargmann transform]
Let $ \vec{\varphi}= ( \varphi_0,..., \varphi_N)$ be a vector-valued function in $L^{2}(\mathbb{R}, \mathbb{H}^{N+1})$. We define $ \mathfrak{B}_{\mathcal{C}}:L^{2}(\mathbb{R}, \mathbb{H}^{N+1}) \to  \mathfrak{A}_{N}(\mathbb{H})$ as
\begin{equation}
\mathfrak{B}_{\mathcal{C}} (\vec{\varphi})(q):= \sum_{n=0}^{N} \mathcal{C}_{N+1} \circ B^{n+1} \varphi_n(q).
\end{equation}
\end{defn}
By using Theorem \ref{CB1} we can write an expression of the $ \mathcal{C}$- full poly Fueter Bargmann transform

\begin{thm}
The $ \mathcal{C}$- full poly Fueter Bargmann transform can be written by using the following commutative diagram
$$\mathfrak{B}_{\mathcal{C}}: \xymatrix{
L^2(\mathbb{R}, \mathbb{H}^{N+1}) \ar[r] \ar[d]_{\mathfrak{B}} & \mathfrak{A}_{N}(\Hq)  \\ \mathcal{\widetilde{F}}_{Slice}^{N+1}(\mathbb{H}) \ar[r]_{Id} & \mathcal{SP}_{N+1}(\Hq) \ar[u]_{ \mathcal{C}_{N+1}}
}$$
with $\mathfrak{B}$ is the QFP Bargmann transform. More precisely, for any vector-valued function $ \vec{\varphi}= ( \varphi_0,..., \varphi_N)$ in $L^{2}(\mathbb{R}, \mathbb{H}^{N+1})$ and $q\in\mathbb{H}$ we have $$ \mathfrak{B}_{\mathcal{C}}(\vec{\varphi})(q)=\int_{\mathbb{R}}\Phi^N(q,x)\varphi_n(x)dx, $$
	where
$$\Phi^N(q,x):=-2^{\frac{3}{4}}\sum_{n=0}^N\sum_{l=0}^{\infty}\sum_{j=0}^{n}\frac{(-1)^j\sqrt{n!}(2\pi)^{\frac{n}{2}}h_{l+n+2}^{2\pi}(x)\mathcal{M}_{n-j,l+n-j}(q,\bar{q})}{2^{\frac{l+n}{2}}(2\pi)^j j!(l+n-j)!(n-j)!}.
$$
\end{thm}

In order to show an unitary and isometric property of the $ \mathcal{C}$- full poly Fueter Bargmann transform we need the following
\begin{defn}
Let us consider the following subspace of $L^{2}(\mathbb{R}, \mathbb{H}^{N+1})$
$$ \mathcal{H}_N:= \bigoplus_{n=0}^N  \bigoplus_{k=n+2}^\infty \{ \psi_{k,n}^{2 \pi} \alpha, \, \alpha \in \mathbb{H}\},$$
where $ \psi_{k,n}^{2 \pi}$ are the components of a vector-valued function.
\end{defn}
\begin{prop}
Let us assume $ \vec{\varphi}, \vec{\psi} \in \mathcal{H}_N$. Then we have
$$ \langle \mathfrak{B}_{\mathcal{C}} \vec{\varphi}, \mathfrak{B}_{\mathcal{C}} \vec{\psi} \rangle_{\mathfrak{A}_{N}(\mathbb{H})}=8 \langle \vec{\varphi}, \vec{\psi} \rangle_{\mathcal{H}_N}.$$
In particular
$$ \| \mathfrak{B}_{\mathcal{C}} \vec{\varphi}\|_{\mathfrak{A}_{N}(\mathbb{H})}=2 \sqrt{2} \| \vec{\varphi}\|_{\mathcal{H}_N}.$$
\end{prop}
\subsection{The map $ \tau_{n+1}$ applied to the QFP Bargmann transform}
Firstly, we study the range of the polyanalytic Fueter mapping $ \tau_{n+1}$ on the quaternionic polyanalytic Fock space,
$$ \mathfrak{\tilde{A}}_{N}(\mathbb{H})= \bigoplus_{n=0}^N \mathcal{\tilde{A}}_{n+1}(\mathbb{H}),$$
where the space $\mathcal{\tilde{A}}_{n+1}(\mathbb{H})$ is defined in \eqref{for2}.
This leads to the following result that extends Theorem \ref{thm24}
\begin{thm}
It holds that
$$ \mathfrak{\tilde{A}}_{N}(\mathbb{H})= \left\{ \sum_{n=0}^N \sum_{h=0}^\infty \mathcal{Q}_{h+n}(q, \bar{q}) \beta_{n,h},   \, \, \sum_{n=0}^N\sum_{h=0}^\infty \frac{(h+n)! | \beta_{n,h}|^2}{4^n n!(h+n+1)(h+n+2)(2 \pi)^{h+2n}}< \infty \right\},$$
where $\{\beta_{h}\}_{h\geq 0} \subset \mathbb{H}$.
\end{thm}

The counterpart of Definition \ref{ad1} in this full polyanalytic setting can be presented as follows

\begin{defn}
Let us consider $ f(q)=\sum_{n=0}^N\sum_{h=0}^\infty \mathcal{Q}_{h+n}(q, \bar{q}) \beta_{n,h}$ and 
\newline
$g(q)= \sum_{n=0}^N\sum_{h=0}^\infty \mathcal{Q}_{h+n}(q, \bar{q}) \gamma_{n,h}$. We define the inner product as
$$ \langle f,g \rangle_{\mathfrak{\tilde{A}}_{N}(\mathbb{H})}= \sum_{n=0}^N \frac{\langle f,g \rangle_{\mathcal{\tilde{A}}_{n+1}(\mathbb{H})}}{4^n (2 \pi)^n n!}=\sum_{n=0}^{N}\sum_{k=0}^\infty \frac{(k+n)!}{4^n n!(k+n+1)(k+n+2)(2 \pi)^{k+2n}} \overline{\gamma_{n,k}} \beta_{n,k}.$$
The associated norm is
$$ \|f\|^2_{\mathfrak{\tilde{A}}_{N}(\mathbb{H})}= \sum_{n=0}^N \frac{\| f\|^2_{\mathcal{\tilde{A}}_{n+1}(\mathbb{H})}}{4^n (2 \pi)^n n!}=\sum_{n=0}^{N}\sum_{k=0}^\infty \frac{(k+n)!}{4^n n!(k+n+1)(k+n+2)(2 \pi)^{k+2n}} |\beta_{n,k}|^2.$$
\end{defn}

\begin{defn}[Full $\tau$- poly Fueter-Bargmann transform]
We define $ \mathfrak{B}_{\tau}: L^{2}(\mathbb{R}, \mathbb{H}^{N+1}) \to \mathfrak{\tilde{A}}_{N}(\mathbb{H})$ as
$$ \mathfrak{B}_{\tau}:=\sum_{n=0}^{N} \tau_{N+1} \circ B^{n+1}.$$
\end{defn}

We present different results in the full polyanalytic setting as follows

\begin{thm}
It is possible to define the full $\mathcal{\tau}$- poly Fueter-Bargmann through this diagram
$$\mathfrak{B}_{\tau}: \xymatrix{
L^2(\mathbb{R}, \mathbb{H}^{N+1}) \ar[r] \ar[d]_{\mathfrak{B}} & \mathfrak{\tilde{A}}_{N}(\mathbb{H})  \\ \mathcal{\widetilde{F}}_{Slice}^{N+1}(\mathbb{H}) \ar[r]_{Id} & \mathcal{SP}_{N+1}(\Hq) \ar[u]_{\tau_{N+1}}
}$$
Precisely, for any $ \vec{\varphi}=(\varphi_0,..., \varphi_N) \in L^{2}(\mathbb{R}, \mathbb{H}^{N+1})$ and $q \in \mathbb{H}$ we have
$$ \mathfrak{B}_{\tau}(\vec{\varphi})(q)= \int_{\mathbb{R}} \Theta_N(q,x) \varphi(x) \, dx,$$
where
\begin{equation}
\Theta_N(q,x)=- 2^{\frac{3}{4}}\sum_{n=0}^{N}  \sum_{h=0}^\infty \frac{\sqrt{n!}  \pi^{\frac{n}{2}}\mathcal{Q}_{k+n}(q, \bar{q}) h_{k+n+2}^{2 \pi}(x)}{2^{\frac{k}{2}}(k+n)!}.
\end{equation}
\end{thm}

\begin{thm}
Let $ \vec{\varphi}$ , $\vec{\psi}$ be vector-valued function in $\mathcal{H}_N$ then
$$ \langle \mathfrak{B}_{\tau} \vec{\varphi}, \mathfrak{B}_{\tau} \vec{\psi}\rangle_{\mathfrak{\tilde{A}}_{N}(\mathbb{H})}=32(2 \pi)^2 \langle \vec{\varphi}, \vec{\psi} \rangle_{ \mathcal{H}_N}.$$
In particular
\begin{equation}
\| \mathfrak{B}_{\tau} \vec{\varphi} \|_{\mathfrak{\tilde{A}}_{N}(\mathbb{H})}= 16  \pi\sqrt{2} \| \vec{\varphi} \|_{ \mathcal{H}_N}.
\end{equation}
\end{thm}

\hspace{4mm}

\noindent
Antonino De Martino,
Dipartimento di Matematica \\ Politecnico di Milano\\
Via Bonardi n.~9\\
20133 Milano\\
Italy

\noindent
\emph{email address}: antonino.demartino@polimi.it\\
\emph{ORCID iD}: 0000-0002-8939-4389

\vspace*{5mm}
\noindent
Kamal Diki,
Schmid College of Science and Technology, \\ Chapman University\\
Orange\\
92866, CA\\
US

\noindent
\emph{email address}: diki@chapman.edu\\
\emph{ORCID iD}: 0000-0002-4359-7535
\end{document}